# CONCENTRATION INEQUALITIES FOR DEPENDENT RANDOM VARIABLES VIA THE MARTINGALE METHOD

By Leonid (Aryeh) Kontorovich[1] and Kavita Ramanan[2]

*Weizmann Institute of Science and Carnegie Mellon University*

The martingale method is used to establish concentration inequalities for a class of dependent random sequences on a countable state space, with the constants in the inequalities expressed in terms of certain mixing coefficients. Along the way, bounds are obtained on martingale differences associated with the random sequences, which may be of independent interest. As applications of the main result, concentration inequalities are also derived for inhomogeneous Markov chains and hidden Markov chains, and an extremal property associated with their martingale difference bounds is established. This work complements and generalizes certain concentration inequalities obtained by Marton and Samson, while also providing different proofs of some known results.

## 1. Introduction.

1.1. *Background.* Concentration of measure is a fairly general phenomenon which, roughly speaking, asserts that a function $\varphi : \Omega \to \mathbb{R}$ with "suitably small" local oscillations defined on a "high-dimensional" probability space $(\Omega, \mathcal{F}, \mathbb{P})$, almost always takes values that are "close" to the average (or median) value of $\varphi$ on $\Omega$. Under various assumptions on the function $\varphi$ and different choices of metrics, this phenomenon has been quite extensively studied in the case when $\mathbb{P}$ is a product measure on a product space $(\Omega, \mathcal{F})$ or, equivalently, when $\varphi$ is a function of a large number of i.i.d. random

Received September 2006; revised July 2007.

[1]Supported in part by NSF Grant DMS-03-23668-0000000965, NSF ITR Grant IIS-0205456.

[2]Supported in part by NSF Grants DMS-04-06191, DMI-0323668-0000000965, DMS-04-05343.

variables (see, e.g., the surveys of Talagrand [29, 30], Ledoux [17], McDiarmid [25] and references therein). Concentration inequalities have found numerous applications in a variety of fields (see, e.g., [6, 25, 28]).

The situation is naturally far more complex for nonproduct measures, where one can trivially construct examples where the concentration property fails. For functions of dependent random variables $(X_i)_{i \in \mathbb{N}}$, the crux of the problem is often to quantify and bound the dependence among the random variables $X_i$, in terms of various types of mixing coefficients. A sufficiently rapid decay of the mixing coefficients often allows one to establish concentration results [22, 23, 27].

A number of techniques have been used to prove measure concentration. Among these are isoperimetric inequalities and the induction method of Talagrand [29, 30], log-Sobolev inequalities developed by Ledoux and others [5, 16, 23, 27], information-theoretic techniques [1, 9, 12, 19, 20, 21, 27], martingale methods based on the Azuma–Hoeffding inequality [2, 8, 11, 25], transportation inequalities (see, e.g., [10, 26]), and Stein's method of exchangeable pairs, recently employed by Chatterjee [7]. The information-theoretic approach has proved quite useful for dealing with nonproduct measures. In a series of papers, Marton [19, 20, 21, 22, 23] successfully used these techniques to establish concentration inequalities for collections of dependent random variables under various assumptions. In this work we adopt a completely different approach, based on the martingale method, to establish concentration bounds for dependent random variables. In the process we establish bounds on certain martingale differences, which may be of independent interest.

In the next subsection we provide a precise description of our main results and discuss their relation to prior work. The subsequent subsections provide an outline of the paper and collect some common notation that we use.

1.2. *Description of main results.* Consider a collection of random variables $(X_i)_{1 \le i \le n}$ taking values in a countable space $\mathcal{S}$. Let $\mathcal{F}$ be the set of all subsets of $\mathcal{S}^n$ and let $\mathbb{P}$ be the probability distribution induced by the finite sequence $X = (X_1, \ldots, X_n)$ on $(\mathcal{S}^n, \mathcal{F})$. Then we can (and will) assume without loss of generality that $X_i, 1 \le i \le n$, are the coordinate projections defined on the probability space $(\mathcal{S}^n, \mathcal{F}, \mathbb{P})$. Given $1 \le i < j \le n$, $x_i^j$ is used to denote the subsequence $(x_i, x_{i+1}, \ldots, x_j)$. Similarly, for $1 \le i < j \le n$, $X_i^j$ represents the random vector $(X_i, \ldots, X_j)$. For further simplicity of notation, $x_1^j$ and $X_1^j$ will be sometimes written simply as $x^j$ and $X^j$, respectively. Let $\mathcal{S}^n$ be equipped with the Hamming metric $d \colon \mathcal{S}^n \times \mathcal{S}^n \to [0, \infty)$, defined by

$$d(x, y) \doteq \sum_{i=1}^{n} \mathbb{1}_{\{x_i \ne y_i\}}.$$

Also, let $\bar{d}(x,y) \doteq d(x,y)/n$ denote the normalized Hamming metric on $\mathcal{S}^n$. In addition, let $\mathbb{E}$ denote expectation with respect to $\mathbb{P}$. Given a function $\varphi \colon \mathcal{S}^n \to \mathbb{R}$, we will use the shorthand notation $\mathbb{P}\{|\varphi - \mathbb{E}\varphi| \geq t\}$ instead of $\mathbb{P}\{|\varphi(X) - \mathbb{E}[\varphi(X)]| \geq t\}$. Also, given two random variables $Y$ and $Z$, $\mathcal{L}(Z \mid Y = y)$ denotes the conditional distribution of $Z$ given $Y = y$.

Our main result is a concentration inequality on the metric probability space $(\mathcal{S}^n, d, \mathbb{P})$, which is expressed in terms of the following mixing coefficients. For $1 \leq i < j \leq n$, define

$$(1.1) \qquad \bar{\eta}_{ij} \doteq \sup_{\substack{y^{i-1} \in \mathcal{S}^{i-1}, w, \hat{w} \in \mathcal{S} \\ \mathbb{P}(X^i = Y^{i-1}w) > 0, \mathbb{P}(X^i = Y^{i-1}\hat{w}) > 0}} \eta_{ij}(y^{i-1}, w, \hat{w}),$$

where, for $y^{i-1} \in \mathcal{S}^{i-1}$ and $w, \hat{w} \in \mathcal{S}$,

$$(1.2) \quad \eta_{ij}(y^{i-1}, w, \hat{w}) \doteq \|\mathcal{L}(X_j^n \mid X^i = y^{i-1}w) - \mathcal{L}(X_j^n \mid X^i = y^{i-1}\hat{w})\|_{\mathrm{TV}}$$

and $\|Q - R\|_{\mathrm{TV}}$ denotes the total variation distance between the probability measures $Q$ and $R$ [see (1.17) for a precise definition]. Moreover, let $\Delta_n$ be the $n \times n$ upper triangular matrix defined by

$$(\Delta_n)_{ij} = \begin{cases} 1, & \text{if } i = j, \\ \bar{\eta}_{ij}, & \text{if } i < j, \\ 0, & \text{otherwise.} \end{cases}$$

Observe that the (usual $\ell_\infty$) operator norm of the matrix $\Delta_n$ is given explicitly by

$$(1.3) \qquad \|\Delta_n\|_\infty = \max_{1 \leq i \leq n} H_{n,i},$$

where, for $1 \leq i \leq n-1$,

$$(1.4) \qquad H_{n,i} \doteq (1 + \bar{\eta}_{i,i+1} + \cdots + \bar{\eta}_{i,n})$$

and $H_{n,n} \doteq 1$.

We can now state the concentration result.

THEOREM 1.1. *Suppose $\mathcal{S}$ is a countable space, $\mathcal{F}$ is the set of all subsets of $\mathcal{S}^n$, $\mathbb{P}$ is a probability measure on $(\mathcal{S}^n, \mathcal{F})$ and $\varphi \colon \mathcal{S}^n \to \mathbb{R}$ is a c-Lipschitz function (with respect to the Hamming metric) on $\mathcal{S}^n$ for some $c > 0$. Then for any $t > 0$,*

$$(1.5) \qquad \mathbb{P}\{|\varphi - \mathbb{E}\varphi| \geq t\} \leq 2 \exp\left(-\frac{t^2}{2nc^2\|\Delta_n\|_\infty^2}\right).$$

Theorem 1.1 follows from Theorem 2.1 and Remark 2.1. For the particular case when $(X_1, \ldots, X_n)$ is a (possibly inhomogeneous) Markov chain, the bound in Theorem 1.1 simplifies further. More precisely, given any initial probability distribution $p_0(\cdot)$ and stochastic transition kernels $p_i(\cdot \mid \cdot)$, $1 \leq i \leq n-1$, let the probability measure $\mathbb{P}$ on $\mathcal{S}^n$ be defined by

$$(1.6) \qquad \mathbb{P}\{(X_1, \ldots, X_i) = x\} = p_0(x_1) \prod_{j=1}^{i-1} p_j(x_{j+1} \mid x_j)$$

for every $1 \leq i \leq n$ and every $x = (x_1, \ldots, x_i) \in \mathcal{S}^i$. Moreover, let $\theta_i$ be the $i$th contraction coefficient of the Markov chain:

$$(1.7) \qquad \theta_i \doteq \sup_{x', x'' \in \mathcal{S}} \|p_i(\cdot \mid x') - p_i(\cdot \mid x'')\|_{\mathrm{TV}}$$

for $1 \leq i \leq n-1$, and set

$$(1.8) \qquad M_n \doteq \max_{1 \leq i \leq n-1} (1 + \theta_i + \theta_i \theta_{i+1} + \cdots + \theta_i \cdots \theta_{n-1}).$$

Then we have the following result.

THEOREM 1.2. *Suppose $\mathbb{P}$ is the Markov measure on $\mathcal{S}^n$ described in (1.6), and $\varphi : \mathcal{S}^n \to \mathbb{R}$ is a c-Lipschitz function with respect to the Hamming metric on $\mathcal{S}^n$ for some $c > 0$. Then for any $t > 0$,*

$$(1.9) \qquad \mathbb{P}\{|\varphi - \mathbb{E}\varphi| \geq t\} \leq 2 \exp\left(-\frac{t^2}{2nc^2 M_n^2}\right),$$

*where $M_n$ is given by (1.8). In addition, if $\varphi$ is c-Lipschitz with respect to the normalized Hamming metric, then*

$$(1.10) \qquad \mathbb{P}\{|\varphi - \mathbb{E}\varphi| \geq t\} \leq 2 \exp\left(-\frac{nt^2}{2c^2 M_n^2}\right).$$

*In particular, when*

$$(1.11) \qquad M \doteq \sup_n \frac{M_n^2}{n} < \infty,$$

*the concentration bound (1.10) is dimension-independent.*

Theorem 1.2 follows from Theorem 1.1 and the observation (proved in Lemma 7.1) that $\|\Delta_n\|_\infty \leq M_n$ for any Markov measure $\mathbb{P}$. In the special case when $\mathbb{P}$ is a uniformly contracting Markov measure, satisfying $\theta_i \leq \theta < 1$ for $1 \leq i \leq n-1$, we have $M_n \leq 1/(1-\theta)$ for every $n$ and the dimension-independence condition (1.11) of Theorem 1.2 also holds trivially. The constants in the exponent of the upper bounds in (1.5) and (1.9) are not sharp—indeed, for the independent case (i.e., when $\mathbb{P}$ is a product measure) it is known that a sharper bound can be obtained by replacing $n/2$ by $2n$ (see [24] and [29]).

Another application of our technique yields a novel concentration inequality for hidden Markov chains. Rather than state the (natural but somewhat long-winded) definitions here, we defer them to Section 7. The key result (Theorem 7.1) is that the mixing coefficients of a hidden Markov chain can be entirely controlled by the contraction coefficients of the underlying Markov chain, and so the concentration bound has the same form as (1.10). The latter result is at least somewhat surprising, as this relationship fails for arbitrary hidden-observed process pairs.

For the purpose of comparison with prior results, it will be useful to derive some simple consequences of Theorem 1.1. If $\varphi$ is a 1-Lipschitz function with respect to the normalized Hamming metric, then it is a $1/n$-Lipschitz function with respect to the regular Hamming metric, and so it follows from Theorem 1.1 that

$$(1.12) \qquad \mathbb{P}\{|\varphi - \mathbb{E}\varphi| \ge t\} \le \alpha(t),$$

where, for $t > 0$,

$$(1.13) \qquad \alpha(t) \doteq 2 \exp\left(-\frac{nt^2}{2\|\Delta_n\|_\infty^2}\right).$$

In turn, this immediately implies the following concentration inequality for any median $m_\varphi$ of $\varphi$: for $t > \|\Delta_n\|_\infty \sqrt{(2\ln 4)/n}$,

$$(1.14) \qquad \mathbb{P}\{|\varphi(X) - m_\varphi| \ge t\} \le 2 \exp\left(-\frac{n}{2}\left(\frac{t}{\|\Delta_n\|_\infty} - \sqrt{\frac{2\ln 4}{n}}\right)^2\right).$$

Indeed, this is a direct consequence of the fact (stated, as Proposition 1.8 of [17]) that if (1.12) holds, then the left-hand side of (1.14) is bounded by $\alpha(t - t_0)$, where $t_0 \doteq \alpha^{-1}(1/2)$. In addition (see, e.g., Proposition 3 of [17]), this also shows that $\alpha_{\mathbb{P}}(\cdot) = \alpha(\cdot - t_0)/2$ acts as a concentration function for the measure $\mathbb{P}$ on $\mathcal{S}^n$ equipped with the normalized Hamming metric: in other words, for any set $A \subset \mathcal{S}^n$ with $\mathbb{P}(A) \ge 1/2$ and $t > t_0$,

$$(1.15) \qquad \mathbb{P}(A_t) \ge 1 - \tfrac{1}{2}\alpha(t - t_0),$$

where $A_t = \{y \in \mathcal{S}^n : \bar{d}(y, x) \le t \text{ for some } x \in A\}$ is the $t$-fattening of $A$ with respect to the normalized Hamming metric.

The above theorems complement the results of Marton [20, 21, 22] and Samson [27]. Theorem 1 of Marton [22] (combined with Lemma 1 of [21] and the comment after Proposition 4 of [20]) shows that when $\mathcal{S}$ is a complete, separable metric space, equipped with the normalized Hamming distance $\bar{d}$, for any Lipschitz function with $\|\varphi\|_{\text{Lip}} \le 1$, the relation (1.15) holds with $t_0 = C\sqrt{(\ln 2)/2n}$ and

$$\alpha(t) \doteq 2 \exp\left(-\frac{2nt^2}{C^2}\right),$$

where

$$C \doteq \max_{1 \le i \le n} \sup_{y \in \mathcal{S}^{i-1}, \bar{w}, \hat{w} \in \mathcal{S}} \inf_{\pi \in \mathcal{M}_i(y, \bar{w}, \hat{w})} \mathbb{E}_\pi [\bar{d}(\bar{X}^n, \hat{X}^n)],$$

with $\mathcal{M}_i(y, \bar{w}, \hat{w})$ being the set of probability measures $\pi = \mathcal{L}(\bar{X}^n, \hat{X}^n)$ on $\mathcal{S}^n \times \mathcal{S}^n$, whose marginals are $\mathcal{L}(X^n \mid X^i = y\bar{w})$ and $\mathcal{L}(X^n \mid X^i = y\hat{w})$, respectively. Moreover, concentration inequalities around the median that are qualitatively similar to the one obtained in Theorem 1.2 were obtained for strictly contracting Markov chains in [20] (see Proposition 1) and for a class of stationary Markov processes in [21] (see Proposition 4'). On the other hand, our result in Theorem 1.2 is applicable to a broader class of Markov chains, which could be nonstationary and not necessarily uniformly contracting.

The mixing coefficients $\bar{\eta}_{ij}$ defined in (1.1) also arise in the work of Samson, who derived concentration bounds for dependent random variables, but in a different space with a different metric, and for a more restrictive class of functions. Specifically, for the case when $\mathcal{S} = [0, 1]$, equipped with the Euclidean metric, it was shown in Samson [27] that if a function $\varphi : \mathcal{S}^n \to \mathbb{R}$ is convex with Lipschitz constant $\|\varphi\|_{\mathrm{Lip}} \le 1$, then

$$\mathbb{P}\{|\varphi - \mathbb{E}\varphi| \ge t\} \le 2 \exp\left(-\frac{t^2}{2\|\Gamma\|_2^2}\right),$$

where $\|\Gamma\|_2$ is the usual $\ell_2$ operator norm of the upper-triangular $n \times n$ matrix $\Gamma$ of the form

$$\Gamma_{ij} \doteq \begin{cases} 1, & \text{if } j = i, \\ (\bar{\eta}_{ij})^{1/2}, & \text{if } i < j, \\ 0, & \text{otherwise.} \end{cases}$$

The results of both Marton and Samson cited above were obtained using a combination of information-theoretic and coupling techniques, as well as the duality method of [4]. In contrast, in this paper we adopt a completely different approach, based on the martingale method (described in Section 2) and a linear algebraic perspective, thus also providing alternative proofs of some known results.

The concentration inequality in Theorem 1.1 was obtained almost contemporaneously with the publication of [8], whose coupling matrix $D^\sigma$ is a close analogue of our $\Delta_n$. We derive essentially the same martingale difference bound as Chazottes et al. by a rather different method—they employ a coupling argument while we rely on the linear programming inequality in Theorem 4.1. The latter is proved in greater generality (for weighted Hamming metrics), and in a much simpler way, in Kontorovich's Ph.D. thesis [15].

1.3. *Outline of paper.* In Section 1.4 we summarize some basic notation used throughout the paper. In Section 2, we set up the basic machinery for applying the martingale method, or equivalently the method of bounded differences, to our problem. In particular, as stated precisely in Theorem 2.1, this reduces the proof of Theorem 1.1 to showing that certain martingale differences $V_i(\varphi; y)$, $1 \le i \le n-1$, associated with the measure $\mathbb{P}$ on $\mathcal{S}^n$ are uniformly bounded by the $\ell_\infty$ operator norm of the matrix $\Delta_n$ of mixing coefficients. Sections 3–5 are devoted to establishing these bounds when $\mathcal{S}$ is finite. The proof uses linear algebraic techniques and establishes a functional inequality that may be of independent interest. Section 6 then uses an approximation argument to extend the bounds to the case when $\mathcal{S}$ is countable. As applications of the main result, Section 7.1 considers the case when $\mathbb{P}$ is a (possibly inhomogeneous) Markov measure, and Section 7.2 performs a similar calculation for measures induced by hidden Markov chains. Specifically, Lemma 7.1 establishes the bound $\|\Delta_n\|_\infty \le M_n$, which then allows Theorem 1.2 to be deduced immediately from Theorem 1.1. Finally, in Section 7.3 we describe a class of extremal functions for martingale differences associated with Markov measures. Some lemmas not central to the paper are collected in the Appendix.

1.4. *Notation and definitions.* In addition to the notation introduced in Section 1.2, we shall use the following common notation throughout the paper. Given a finite or countable set $\mathcal{S}$ and finite sequences $x \in \mathcal{S}^k$ and $y \in \mathcal{S}^\ell$, we use either $xy \in \mathcal{S}^{k+\ell}$ or $[x\,y] \in \mathcal{S}^{k+\ell}$ to denote concatenation of the two sequences. The space $\mathcal{S}^0$ represents the null string. Also, we will use the shorthand notation $\sum_{x_i^j}$ to mean $\sum_{x_i^j \in \mathcal{S}^{j-i+1}}$. The random variables $X = (X_1, \dots, X_n)$ defined on $(\mathcal{S}^n, \mathcal{F}, \mathbb{P})$ will always represent coordinate projections: $X_i(x) = x_i$. Therefore for conciseness, we will sometimes simply write $\mathbb{P}(x)$ and $\mathbb{P}(x_j^n \mid x^i)$ to denote $\mathbb{P}\{X = x\}$ and $\mathbb{P}\{X_j^n = x_j^n \mid X^i = x^i\}$, respectively. As in (1.1), we will always assume (often without explicit mention) that terms involving conditional probabilities are restricted to elements for which the probability of the conditioned event is strictly positive.

The indicator variable $\mathbb{1}_{\{\cdot\}}$ takes on the value 1 if the predicate in the bracket $\{\cdot\}$ is true, and 0 otherwise. The sign function is defined by $\operatorname{sgn}(z) = \mathbb{1}_{\{z > 0\}} - \mathbb{1}_{\{z \le 0\}}$ and the positive function is defined by $(z)_+ = \max(z, 0) = z \mathbb{1}_{\{z > 0\}}$. We use the standard convention that $\sum_{z \in A} z = 0$ and $\prod_{z \in A} z = 1$ whenever $A$ is empty ($A = \varnothing$).

Throughout the paper, $K_n$ denotes the space of all functions $\kappa \colon \mathcal{S}^n \to \mathbb{R}$ (for finite $\mathcal{S}$) and $\Phi_n \subset K_n$ the subset of 1-Lipschitz functions $\varphi \colon \mathcal{S}^n \to [0, n]$.

For a discrete, signed measure space $(\mathcal{X}, \mathcal{B}, \nu)$, recall that the $\ell_1$ norm is given by

$$(1.16) \qquad \|\nu\|_1 = \sum_{x \in \mathcal{X}} |\nu(x)|.$$

Given two probability measures $\nu_1$, $\nu_2$ on the measurable space $(\mathcal{X}, \mathcal{B})$, we define the total variation distance between the two measures as follows:

$$(1.17) \qquad \|\nu_1 - \nu_2\|_{TV} = \sup_{A \in \mathcal{B}} |\nu_1(A) - \nu_2(A)|$$

(which, by common convention, is equal to *half* the total variation of the signed measure $\nu_1 - \nu_2$). It is easy to see that in this case,

$$(1.18) \qquad \|\nu_1 - \nu_2\|_{TV} = \tfrac{1}{2} \|\nu_1 - \nu_2\|_1 = \sum_{x \in \mathcal{X}} (\nu_1(x) - \nu_2(x))_+.$$

## 2. Method of bounded martingale differences.

Since our proof relies on the so-called martingale method of establishing concentration inequalities, here we briefly introduce the method (see [24, 25] or Section 4.1 of [17] for more thorough treatments). Let $X = (X_i)_{1 \leq i \leq n}$ be a collection of random variables defined on a probability space $(\Omega, \mathcal{F}, \mathbb{P})$, taking values in a space $\mathcal{S}$. Then given any filtration of sub-$\sigma$-algebras,

$$\{\varnothing, \Omega\} = \mathcal{F}_0 \subset \mathcal{F}_1 \subset \cdots \subset \mathcal{F}_n = \mathcal{F},$$

and a function $\varphi: \mathcal{S}^n \to \mathbb{R}$, define the associated martingale differences by

$$(2.1) \qquad V_i(\varphi) = \mathbb{E}[\varphi(X) \mid \mathcal{F}_i] - \mathbb{E}[\varphi(X) \mid \mathcal{F}_{i-1}]$$

for $i = 1, \ldots, n$. It is a classical result, going back to Hoeffding [11] and Azuma [2], that

$$(2.2) \qquad \mathbb{P}\{|\varphi - \mathbb{E}\varphi| \geq r\} \leq 2 \exp(-r^2/2D^2)$$

for any $D$ such that $D^2 \geq \sum_{i=1}^{n} \|V_i(\varphi)\|_\infty^2$.

In the setting of this paper, we have $(\Omega, \mathcal{F}, \mathbb{P}) = (\mathcal{S}^n, \mathcal{F}, \mathbb{P})$, where $\mathcal{S}$ is a countable set, $\mathcal{F}$ is the set of all subsets of $\mathcal{S}^n$ and $X = (X_i)_{1 \leq i \leq n}$ is the collection of coordinate projections. For $i = 1, \ldots, n$, we set $\mathcal{F}_0 \doteq \{\varnothing, \mathcal{S}^n\}$, $\mathcal{F}_n \doteq \mathcal{F}$ and for $1 \leq i \leq n-1$, let $\mathcal{F}_i$ be the $\sigma$-algebra generated by $X^i = (X_1, \ldots, X_i)$. Given any function $\varphi$ on $\mathcal{S}^n$, for $1 \leq i \leq n$, define the martingale differences $V_i(\varphi)$ in the standard way, by (2.1).

The following theorem shows that when $\varphi$ is Lipschitz, these martingale differences can be bounded in terms of the mixing coefficients defined in Section 1.2. The proof of the theorem is given in Section 6.

THEOREM 2.1. *If $\mathcal{S}$ is a countable set and $\mathbb{P}$ is a probability measure on $(\mathcal{S}^n, \mathcal{F})$ such that $\min_{i=1,\ldots,n} \inf_{y^i \in \mathcal{S}^i : \mathbb{P}(X^i = y^i) > 0} \mathbb{P}(X^i = y^i) > 0$, then for any 1-Lipschitz function $\varphi$ on $\mathcal{S}^n$, we have for $1 \leq i \leq n$,*

$$(2.3) \qquad \|V_i(\varphi)\|_\infty \leq H_{n,i} = 1 + \sum_{j=i+1}^{n} \bar{\eta}_{ij},$$

*where $\{V_i(\varphi)\}$ are the martingale differences defined in (2.1) and the coefficients $\bar{\eta}_{ij}$ and $H_{n,i}$ are as defined in (1.1) and (1.4), respectively.*

REMARK 2.1. Since $\|V_i(\cdot)\|_\infty$ is homogeneous in the sense that $\|V_i(a\varphi)\|_\infty = a\|V_i(\varphi)\|_\infty$ for $a > 0$, and $a^{-1}\varphi$ is 1-Lipschitz whenever $\varphi$ is $a$-Lipschitz, Theorem 2.1 implies that for $1 \leq i \leq n$,

$$\|V_i(\varphi)\|_\infty \leq cH_{n,i}$$

for any $c$-Lipschitz $\varphi$. Along with the relation (1.3), this shows that

$$\sum_{i=1}^n \|V_i(\varphi)\|_\infty^2 \leq n \max_{1 \leq i \leq n} \|V_i(\varphi)\|_\infty^2 \leq nc^2 \max_{1 \leq i \leq n} H_{n,i}^2 = nc^2\|\Delta_n\|_\infty^2.$$

When combined with Azuma's inequality (2.2), this shows that Theorem 2.1 implies Theorem 1.1.

REMARK 2.2. The quantity $V_i(\cdot)$ is translation-invariant in the sense that $V_i(\varphi + a) = V_i(\varphi)$ for every $a \in \mathbb{R}$. Since the length of the range of any 1-Lipschitz function on $\mathcal{S}^n$ is equal to $n$, the Hamming diameter of $\mathcal{S}^n$, for the proof of Theorem 2.1 there is no loss of generality in assuming that the range of $\varphi$ lies in $[0, n]$.

**3. A linear programming bound for martingale differences.** In the next three sections, we prove Theorem 2.1 under the assumption that $\mathcal{S}$ is finite. We start by obtaining a slightly more tractable form for the martingale difference.

LEMMA 3.1. *Given a probability measure $\mathbb{P}$ on $(\mathcal{S}^n, \mathcal{F})$ and any function $\varphi : \mathcal{S}^n \to \mathbb{R}$, let the martingale differences $\{V_i(\varphi), 1 \leq i \leq n\}$ be defined as in* (2.1)*. Then, for $1 \leq i \leq n$,*

$$\|V_i(\varphi)\|_\infty \leq \max_{\substack{y^{i-1} \in \mathcal{S}^{i-1}, w, \hat{w} \in \mathcal{S} \\ \mathbb{P}(X^i = Y^{i-1}w) > 0, \mathbb{P}(X^i = Y^{i-1}\hat{w}) > 0}} |\hat{V}_i(\varphi; y^{i-1}, w, \hat{w})|,$$

*where, for $y^{i-1} \in \mathcal{S}^{i-1}$ and $w, \hat{w} \in \mathcal{S}$,*

$$(3.1) \quad \hat{V}_i(\varphi; y^{i-1}, w, \hat{w}) \doteq \mathbb{E}[\varphi(X) \mid X^i = y^{i-1}w] - \mathbb{E}[\varphi(X) \mid X^i = y^{i-1}\hat{w}].$$

PROOF. Since $V_i(\varphi)$ is $\mathcal{F}_i$-measurable and $\mathcal{F}_i = \sigma(X^i)$, it follows immediately that

$$(3.2) \qquad\qquad \|V_i(\varphi)\|_\infty = \max_{z^i \in \mathcal{S}^i} |V_i(\varphi; z^i)|,$$

where for $1 \leq i \leq n$ and $z^i \in \mathcal{S}^i$, we define

$$(3.3) \qquad V_i(\varphi; z^i) = \mathbb{E}[\varphi(X) \mid X^i = z^i] - \mathbb{E}[\varphi(X) \mid X^{i-1} = z^{i-1}].$$

Expanding $V_i(\varphi; z^i)$, we obtain

$$V_i(\varphi; z^i)$$
$$= \mathbb{E}[\varphi(X) \mid X^i = z^i] - \sum_{\hat{w} \in \mathcal{S}} \mathbb{E}[\varphi(X) \mid X^i = z^{i-1}\hat{w}]\mathbb{P}(z^{i-1}\hat{w} \mid z^{i-1})$$
$$= \sum_{\hat{w} \in \mathcal{S}} \mathbb{P}(z^{i-1}\hat{w} \mid z^{i-1})(\mathbb{E}[\varphi(X) \mid X^i = z^i] - \mathbb{E}[\varphi(X) \mid X^i = z^{i-1}\hat{w}])$$
$$= \sum_{\hat{w} \in \mathcal{S}} \mathbb{P}(z^{i-1}\hat{w} \mid z^{i-1})\hat{V}_i(\varphi; z^{i-1}, z_i, \hat{w}),$$

where the second equality uses the fact that $\sum_{\hat{w} \in \mathcal{S}} \mathbb{P}(z^{i-1}\hat{w} \mid z^{i-1}) = 1$ with, as usual, $\mathbb{P}(z^{i-1}\hat{w} \mid z^{i-1})$ representing $\mathbb{P}\{X^i = z^{i-1}\hat{w} \mid X^{i-1} = z^{i-1}\}$. In turn, since $0 \le \mathbb{P}(z^{i-1}\hat{w} \mid z^{i-1}) \le 1$, the last display implies (via Jensen's inequality) that for every $z^i \in \mathcal{S}^i$,

$$|V_i(\varphi; z^i)| \le \sum_{\hat{w} \in \mathcal{S}} \mathbb{P}(z^{i-1}\hat{w} \mid z^{i-1})|\hat{V}_i(\varphi; z^{i-1}, z_i, \hat{w})| \le \max_{\hat{w} \in \mathcal{S}} |\hat{V}_i(\varphi; z^{i-1}, z_i, \hat{w})|.$$

Taking the maximum of both sides over $z^i \in \mathcal{S}^i$ and invoking (3.2), the desired inequality is obtained. $\square$

For $n \in \mathbb{N}$, define the finite-dimensional vector space

$$(3.4) \qquad K_n \doteq \{\kappa : \mathcal{S}^n \to \mathbb{R}\}$$

which becomes a Euclidean space when endowed with the inner product

$$(3.5) \qquad \langle \kappa, \lambda \rangle \doteq \sum_{x \in \mathcal{S}^n} \kappa(x)\lambda(x).$$

Also, let $K_0$ be the collection of scalars.

Now, note that for $y^{n-1} \in \mathcal{S}^{n-1}$ and $w, \hat{w} \in \mathcal{S}$,

$$\hat{V}_n(\varphi; y^{n-1}, w, \hat{w}) = \varphi(y^{n-1}w) - \varphi(y^{n-1}\hat{w}),$$

and thus for all 1-Lipschitz functions $\varphi$, the bound

$$(3.6) \qquad |\hat{V}_n(\varphi; y^{n-1}, w, \hat{w})| \le 1 = H_{n,n}$$

holds immediately. On the other hand, given $1 \le i \le n - 1$, $y^{i-1} \in \mathcal{S}^{i-1}$ and $w, \hat{w} \in \mathcal{S}$, the map $\varphi \mapsto \hat{V}_i(\varphi; y^{i-1}, w, \hat{w})$ defined in (3.1) is clearly a linear functional on $K_n$. It therefore admits a representation as an inner product with some element $\kappa = \kappa[y^{i-1}, w, \hat{w}] \in K_n$ (where the notation $[\cdot]$ is used to emphasize the dependence of the function $\kappa \in K_n$ on $y^{i-1}, w$ and $\hat{w}$):

$$(3.7) \qquad \hat{V}_i(\varphi; y^{i-1}, w, \hat{w}) = \langle \kappa, \varphi \rangle.$$

Indeed, expanding (3.1), it is easy to see that the precise form of $\kappa = \kappa[y^{i-1}, w, \hat{w}]$ is given by

$$(3.8) \quad \kappa(x) = \mathbb{1}_{\{x^i = y^{i-1}w\}} \mathbb{P}(x^n_{i+1} \mid y^{i-1}w) - \mathbb{1}_{\{x^i = y^{i-1}\hat{w}\}} \mathbb{P}(x^n_{i+1} \mid y^{i-1}\hat{w})$$

for $x \in \mathcal{S}^n$.

Define $\Phi_n \subset K_n$ to be the subset of 1-Lipschitz functions (with respect to the Hamming metric) with range in $[0, n]$. As observed in Remark 2.2, in order to prove Theorem 2.1, it suffices to establish the martingale difference bounds (2.3) just for $\varphi \in \Phi_n$. In light of Lemma 3.1 and the representation (3.7), it is therefore natural to study the quantity

$$(3.9) \quad \|\kappa\|_\Phi = \max_{\varphi \in \Phi_n} |\langle \kappa, \varphi \rangle|$$

for $\kappa \in K_n$. The notation used reflects the fact that $\|\cdot\|_\Phi$ defines a norm on $K_n$. Since we will not be appealing to any properties of norms, we relegate the proof of this fact to Appendix A.1.

REMARK 3.2. The title of this section is motivated by the fact that the optimization problem $\max_{\varphi \in \Phi_n} \langle \kappa, \varphi \rangle$ is a linear program. Indeed, $f(\cdot) = \langle \kappa, \cdot \rangle$ is a linear function and $\Phi_n$ is a finitely generated, compact, convex polytope. We make no use of this simple fact in our proofs and therefore do not prove it, but see Lemma 4.4 for a proof of a closely related claim.

## 4. A bound on the $\Phi$-norm.

In Section 3 we motivated the introduction of the norm $\|\cdot\|_\Phi$ on the space $K_n$. In this section we bound $\|\cdot\|_\Phi$ by another, more tractable, norm, which we call the $\Psi$-norm. In Section 5 we then bound the $\Psi$-norm in terms of the coefficients $H_{n,i}$.

For $n \in \mathbb{N}$, define the *marginal projection operator* $(\cdot)'$, which takes $\kappa \in K_n$ to $\kappa' \in K_{n-1}$ as follows: if $n > 1$, for each $y \in \mathcal{S}^{n-1}$,

$$(4.1) \quad \kappa'(y) \doteq \sum_{x_1 \in \mathcal{S}} \kappa(x_1 y);$$

if $n = 1$, then $\kappa'$ is the scalar

$$\kappa' = \sum_{x_1 \in \mathcal{S}} \kappa(x_1).$$

We define the Positive-Summation-Iterated (Psi) functional $\Psi_n \colon K_n \to \mathbb{R}$ recursively using projections: $\Psi_0(\cdot) = 0$ and for $n \geq 1$,

$$(4.2) \quad \Psi_n(\kappa) \doteq \sum_{x \in \mathcal{S}^n} (\kappa(x))_+ + \Psi_{n-1}(\kappa'),$$

where we recall that $(z)_+ = \max(z, 0)$ is the positive part of $z$. The norm associated with $\Psi_n$ is then defined to be

$$(4.3) \quad \|\kappa\|_\Psi \doteq \max_{s \in \{-1, 1\}} \Psi_n(s\kappa).$$

As in the case of $\|\cdot\|_\Phi$, it is easily verified that $\|\cdot\|_\Psi$ is a valid norm on $K_n$ (see Lemma A.1).

The next theorem is the main result of this section.

THEOREM 4.1.  *For all $n \in \mathbb{N}$ and $\kappa \in K_n$,*

$$\|\kappa\|_\Phi \leq \|\kappa\|_\Psi.$$

The remainder of the section is devoted to proving this theorem. See Kontorovich's Ph.D. thesis [15] for a considerably simpler proof of a more general claim (which covers the weighted Hamming metrics).

REMARK 4.2.  We will assume for simplicity that $z \neq 0$ whenever expressions of the form $\mathrm{sgn}(z)$ or $\mathbb{1}_{\{z>0\}}$ are encountered below. This incurs no loss of generality [the inequalities proved for this special case will hold in general by continuity of $(\cdot)_+$] and affords us a slightly cleaner exposition, obviating the need to check the $z = 0$ case.

First, we need to introduce a bit more notation. For $n \in \mathbb{N}$ and $y \in \mathcal{S}$, define the *$y$-section* operator $(\cdot)_y \colon K_n \to K_{n-1}$ that takes $\kappa$ to $\kappa_y$ by

$$(4.4) \qquad \kappa_y(x) \doteq \kappa(xy)$$

for $x \in \mathcal{S}^{n-1}$. By convention, for $n = 1$ and $x \in \mathcal{S}^0$, $\kappa_y(x)$ is equal to the scalar $\kappa_y = \kappa(y)$.

Note that for any $y \in \mathcal{S}$, the marginal projection and $y$-section operators commute; in other words, for $\kappa \in K_{n+2}$, we have $(\kappa')_y = (\kappa_y)' \in K_n$ and so we can denote this common value simply by $\kappa'_y \in K_n$: for each $z \in \mathcal{S}^n$,

$$(4.5) \qquad \kappa'_y(z) = \sum_{x_1 \in \mathcal{S}} \kappa_y(x_1 z) = \sum_{x_1 \in \mathcal{S}} \kappa(x_1 z y).$$

Moreover, summing both sides of the first equality in (4.5) over $z \in \mathcal{S}^n$, we obtain

$$(4.6) \qquad \sum_{z \in \mathcal{S}^n} \kappa'_y(z) = \sum_{z \in \mathcal{S}^n} \sum_{x_1 \in \mathcal{S}} \kappa_y(x_1 z) = \sum_{x \in \mathcal{S}^{n+1}} \kappa_y(x).$$

We can use $y$-sections to recast the $\Psi_n(\cdot)$ functional in an alternative form:

LEMMA 4.3.  *For all $n \geq 0$ and $\kappa \in K_{n+1}$, we have*

$$(4.7) \qquad \Psi_{n+1}(\kappa) = \sum_{y \in \mathcal{S}} \left[ \Psi_n(\kappa_y) + \left( \sum_{x \in \mathcal{S}^n} \kappa_y(x) \right)_+ \right].$$

Proof.    Suppose $n = 0$. Then for $\kappa \in K_1$, $\kappa_y \in K_0$ is the scalar $\kappa(y)$, $\Psi_0(\kappa_y) = 0$ by definition and $\mathcal{S}^0$ consists of just one element (the null string). Thus the r.h.s. of (4.7) becomes $\sum_{y \in \mathcal{S}} (\kappa(y))_+$, which by definition [see (4.2)] is equal to $\Psi_1(\kappa)$. So the claim holds for $n = 0$.

Now suppose (4.7) holds for $n = 0, \ldots, N$ for some $N \geq 0$. In order to prove the claim for $n = N + 1$, we pick any $\lambda \in K_{N+2}$ and observe that

$$
\begin{aligned}
\sum_{y \in \mathcal{S}} & \left[ \Psi_{N+1}(\lambda_y) + \left( \sum_{x \in \mathcal{S}^{N+1}} \lambda_y(x) \right)_+ \right] \\
(4.8) \quad & = \sum_{y \in \mathcal{S}} \left[ \left( \Psi_N(\lambda'_y) + \sum_{x \in \mathcal{S}^{N+1}} (\lambda_y(x))_+ \right) + \left( \sum_{x \in \mathcal{S}^{N+1}} \lambda_y(x) \right)_+ \right] \\
& = \sum_{y \in \mathcal{S}} \left[ \Psi_N(\lambda'_y) + \left( \sum_{u \in \mathcal{S}^N} \lambda'_y(u) \right)_+ \right] + \sum_{z \in \mathcal{S}^{N+2}} (\lambda(z))_+
\end{aligned}
$$

where the first equality results from the definition of $\Psi_n$ in (4.2) and the second equality uses the trivial identity

$$
\sum_{y \in \mathcal{S}} \sum_{x \in \mathcal{S}^{N+1}} (\lambda_y(x))_+ = \sum_{z \in \mathcal{S}^{N+2}} (\lambda(z))_+
$$

and the relation (4.6) (with $\kappa$ replaced by $\lambda$ and $n$ by $N$).

On the other hand, by the definition given in (4.2) we have

$$
(4.9) \qquad \Psi_{N+2}(\lambda) = \sum_{z \in \mathcal{S}^{N+2}} (\lambda(z))_+ + \Psi_{N+1}(\lambda').
$$

To compare the r.h.s. of (4.8) with the r.h.s. of (4.9), note that the term $\sum_{z \in \mathcal{S}^{N+2}} (\lambda(z))_+$ is common to both sides and, since (4.7) is satisfied with $n = N$ by the inductive hypothesis, the remaining two terms are also equal:

$$
\Psi_{N+1}(\lambda') = \sum_{y \in \mathcal{S}} \left[ \Psi_N(\lambda'_y) + \left( \sum_{u \in \mathcal{S}^N} \lambda'_y(u) \right)_+ \right].
$$

This establishes (4.7) for $n = N + 1$ and the lemma follows by induction.    $\square$

We will need one more definition to facilitate the main proof. Fix a function $\kappa \in K_n$ and consider some properties that any other function $\alpha \in K_n$ might have relative to $\kappa$:

| | |
|---|---|
| (SL1) | for all $x \in \mathcal{S}^n$, <br> $0 \leq \alpha(x) \leq n - \mathbb{1}_{\{\kappa(x) > 0\}}$ |
| (SL1$'$) | for all $x \in \mathcal{S}^n$, <br> $0 \leq \alpha(x) \leq n - \mathbb{1}_{\{\kappa(x) > 0\}} + 1$ |
| (SL2) | for all $x, y \in \mathcal{S}^n$, <br> $\operatorname{sgn} \kappa(x) = \operatorname{sgn} \kappa(y) \implies |\alpha(x) - \alpha(y)| \leq d(x,y)$ |
| (SL3) | for all $x, y \in \mathcal{S}^n$ with $d(x,y) = 1$, <br> $\operatorname{sgn} \kappa(x) > \operatorname{sgn} \kappa(y) \implies \alpha(x) \leq \alpha(y) \leq \alpha(x) + 2$ |

We define $A_n(\kappa)$ to be the set of $\alpha \in K_n$ that satisfy (SL1), (SL2) and (SL3). Similarly, we write $B_n(\kappa)$ to denote the set of all $\beta \in K_n$ that satisfy (SL1$'$), (SL2) and (SL3). (A possible descriptive name for these objects is "sub-Lipschitz polytopes"—hence the letters SL.)

The following is almost an immediate observation.

LEMMA 4.4.   *For $n \in \mathbb{N}$ and $\kappa \in K_n$, the following two properties hold:*

(a)  *$A_n(\kappa)$ and $B_n(\kappa)$ are compact, convex polytopes in $[0, n+1]^{\mathcal{S}^n}$;*
(b)  *for all $y \in \mathcal{S}$,*

$$(4.10) \qquad \kappa\alpha \in A_n(\kappa) \implies \alpha_y \in B_{n-1}(\kappa_y).$$

PROOF.   Property (a) is verified by checking that each of (SL1), (SL1$'$), (SL2) and (SL3) is closed under convex combinations. To verify property (b), fix $y \in \mathcal{S}$ and choose $\alpha \in A_n(\kappa)$. Using the definition $\alpha_y(x) = \alpha(xy)$ and the fact that $d(x,z) = d(xy, zy)$ for $x, z \in \mathcal{S}^{n-1}$, it is straightforward to check that the fact that $\alpha$ satisfies (SL2) [resp., (SL3)] relative to $\kappa$ implies that $\alpha_y$ also satisfies the same property relative to $\kappa_y$. Moreover, since

$$n - \mathbb{1}_{\{\kappa(x) > 0\}} = (n-1) - \mathbb{1}_{\{\kappa(x) > 0\}} + 1,$$

the fact that $\alpha$ satisfies (SL1) relative to $\kappa$ implies that $\alpha_y$ satisfies (SL1$'$) relative to $\kappa_y$. This proves property (b) and hence the lemma.   $\square$

We will also need the following simple fact about $B_n(\kappa)$ [which, in fact, also holds for $A_n(\kappa)$].

LEMMA 4.5.   *For $n \in \mathbb{N}$ and any $\kappa \in K_n$, if $\hat{\beta}$ is an extreme point of $B_n(\kappa)$, then $\hat{\beta}(x)$ is an integer between $0$ and $n+1$ for every $x \in \mathcal{S}^n$ [in other words, $\hat{\beta} \in B_n(\kappa) \cap \mathbb{Z}_+^{\mathcal{S}^n}$].*

PROOF.   Fix $n \in \mathbb{N}$ and $\kappa \in K_n$. We will establish the lemma using an argument by contradiction. Suppose that $B_n(\kappa)$ has an extreme point $\beta$ that takes on a noninteger value for some $x \in \mathcal{S}^n$. Let $E \subseteq \mathcal{S}^n$ be the set of elements of $\mathcal{S}^n$ on which $\beta$ is not an integer:

$$E \doteq \{x \in \mathcal{S}^n : \beta(x) \notin \mathbb{N}\}.$$

Define $\beta_\varepsilon^+$ by

$$\beta_\varepsilon^+(x) \doteq \beta(x) + \varepsilon \mathbb{1}_{\{x \in E\}}, \tag{4.11}$$

and similarly

$$\beta_\varepsilon^-(x) \doteq \beta(x) - \varepsilon \mathbb{1}_{\{x \in E\}}, \tag{4.12}$$

where $\varepsilon \in (0, 1/2)$ is chosen small enough to satisfy

$$\lfloor \beta(x) \rfloor \leq \tilde{\beta}(x) \leq \lceil \beta(x) \rceil, \tag{4.13}$$

for $\tilde{\beta} = \beta_\varepsilon^-$ and $\tilde{\beta} = \beta_\varepsilon^+$ and

$$\beta(x) < \beta(y) \quad \Longrightarrow \quad \beta(x) + \varepsilon < \beta(y) \tag{4.14}$$

for all $x, y \in \mathcal{S}^n$ (such a choice of $\varepsilon$ is feasible because $\mathcal{S}^n$ is a finite set).

We claim that $\beta_\varepsilon^+ \in B_n(\kappa)$. If $x \notin E$, then $\beta_\varepsilon^+(x) = \beta_\varepsilon^-(x) = \beta(x)$. On the other hand, if $x \in E$, then $(\mathtt{SL1'})$ must hold with strict inequalities:

$$0 < \beta(x) < n - \mathbb{1}_{\{\kappa(x)>0\}} + 1.$$

Together with the condition (4.13), this ensures that $\beta_\varepsilon^+$ and $\beta_\varepsilon^-$ also satisfy $(\mathtt{SL1'})$. Similarly, the relation (4.14) can be used to infer that $\beta_\varepsilon^+$ and $\beta_\varepsilon^-$ satisfy $(\mathtt{SL3})$. It only remains to verify $(\mathtt{SL2})$. First observe that $\beta_\varepsilon^+(x) - \beta_\varepsilon^+(y) = \beta(x) - \beta(y)$ whenever $\{x, y\} \subseteq E$, $\{x, y\} \subseteq E^c$. Thus to prove $(\mathtt{SL2})$, by symmetry it suffices to only consider $x, y \in \mathcal{S}^n$ such that $\mathrm{sgn}\,\kappa(x) = \mathrm{sgn}\,\kappa(y)$, $x \in E$ and $y \notin E$. In this case, we have $\beta(x) \neq \beta(y)$ and

$$\beta_\varepsilon^+(x) - \beta_\varepsilon^+(y) = \beta(x) - \beta(y) + \varepsilon. \tag{4.15}$$

If $\beta(x) < \beta(y)$, then (4.14) and the fact that $\beta$ satisfies $(\mathtt{SL2})$ show that

$$-d(x, y) \leq \varepsilon - d(x, y) \leq \beta(x) - \beta(y) + \varepsilon < 0.$$

The last two displays, together, show that then $\beta_\varepsilon^+$ satisfies $(\mathtt{SL2})$. On the other hand, suppose $\beta(x) > \beta(y)$. Since $d(\cdot, \cdot)$ is the Hamming metric and $x \neq y$, $d(x, y)$ is an integer greater than or equal to 1. Moreover, $\beta(y)$ is also an integer since $y \notin E$. Together with the fact that $\beta$ satisfies $(\mathtt{SL2})$ this implies that $\beta(x) < \beta(y) + d(x, y)$. Therefore, since $\mathcal{S}$ is finite, by choosing $\varepsilon > 0$ smaller if necessary, one can assume that $\beta_\varepsilon^+(x) \leq \beta_\varepsilon^+(y) + d(x, y)$. When combined with the elementary inequality $\beta_\varepsilon^+(x) - \beta_\varepsilon^+(y) = \beta(x) + \varepsilon - \beta(y) > 0$, this proves that $\beta_\varepsilon^+$ satisfies $(\mathtt{SL2})$. The argument for $\beta_\varepsilon^-$ is analogous and thus omitted.

However, having proved that $\beta_\varepsilon^-, \beta_\varepsilon^+ \in B_n(\kappa)$, we now see that $\beta = \frac{1}{2}\beta_\varepsilon^+ + \frac{1}{2}\beta_\varepsilon^-$ is a strict convex combination of two elements of $B_n(\kappa)$, which contradicts the fact that it is an extreme point of this set.  □

Let us observe that $A_n(\kappa) \subset B_n(\kappa)$. More importantly, we will utilize the structural relationship between $A_n$ and $B_n$ stated in the next lemma.

LEMMA 4.6.  *Let $\kappa \in K_n$ be given. For any $\beta \in B_n(\kappa)$ there is an $\alpha \in A_n(\kappa)$ such that*

$$(4.16) \qquad \langle \kappa, \beta \rangle \leq \left( \sum_{x \in \mathcal{S}^n} \kappa(x) \right)_+ + \langle \kappa, \alpha \rangle.$$

PROOF.  Fix $\kappa \in K_n$ throughout the proof. Since we are trying to bound $\langle \kappa, \beta \rangle$ for $\beta \in B_n(\kappa)$ and linear functions achieve their maxima on the extreme points of convex sets, and $A_n(\kappa)$ is convex, it suffices to establish (4.16) only for $\beta$ that are extreme points of $B_n(\kappa)$. By Lemma 4.5, this implies that $\beta \in B_n(\kappa) \cap \mathbb{Z}_+^{\mathcal{S}^n}$. If, in addition, $\beta \in A_n(\kappa)$, the statement of the lemma is trivial; so we need only consider

$$(4.17) \qquad \beta \in (B_n(\kappa) \setminus A_n(\kappa)) \cap \mathbb{Z}_+^{\mathcal{S}^n}.$$

For $\beta$ satisfing (4.17), define $\alpha$ by

$$(4.18) \qquad \alpha(x) \doteq (\beta(x) - 1)_+ \qquad \text{for } x \in \mathcal{S}^n.$$

We first claim that $\alpha \in A_n(\kappa)$. To see why this is true, first observe that the fact that $\beta$ satisfies (SL1′) immediately implies that $\alpha$ satisfies (SL1). Moreover, for $x, y \in \mathcal{S}^n$, if $\beta(x) \leq 1, \beta(y) \leq 1$, then $\alpha(x) - \alpha(y) = 0$; if $\beta(x) > 1, \beta(y) > 1$, then $\alpha(x) - \alpha(y) = \beta(x) - \beta(y)$; if $\beta(x) \leq 1, \beta(y) > 1$, then

$$0 \geq \alpha(x) - \alpha(y) = -\beta(y) + 1 \geq \beta(x) - \beta(y),$$

with an analogous relation holding if $\beta(x) > 1, \beta(y) \leq 1$.

Combining the above relations with the fact that $\beta$ satisfies properties (SL2) and (SL3), it is straightforward to see that $\alpha$ also satisfies the same properties. Having shown that $\alpha \in A_n(\kappa)$, we will now show that $\alpha$ satisfies (4.16). Proving (4.16) is equivalent to showing that

$$(4.19) \qquad \langle \kappa, \delta \rangle = \sum_{x \in \mathcal{S}^n} \kappa(x)\delta(x) \leq \left( \sum_{x \in \mathcal{S}^n} \kappa(x) \right)_+$$

where $\delta(x) = \beta(x) - \alpha(x) = \mathbb{1}_{\{\beta(x) \geq 1\}}$. To this end, we claim (and justify below) that for $\beta$ satisfying (4.17) and for $\delta$ defined as above, we have for all $z \in \mathcal{S}^n$

$$(4.20) \qquad \kappa(z) < 0 \quad \Longrightarrow \quad \delta(z) = 1;$$

the relation (4.19) then easily follows from (4.20).

Suppose, to get a contradiction, that (4.20) fails; this means that there exists a $z \in \mathcal{S}^n$ such that

$$(4.21) \qquad \kappa(z) < 0, \qquad \beta(z) = 0.$$

Recall that by assumption, $\beta \notin A_n(\kappa)$; since $\beta \in B_n(\kappa)$, this can only happen if $\beta$ violates (SL1), which can occur in one of only two ways:

(i) $\beta(y) = n + 1$ for some $y \in \mathcal{S}^n$ with $\kappa(y) \leq 0$;
(ii) $\beta(y) = n$ for some $y \in \mathcal{S}^n$ with $\kappa(y) > 0$.

We can use (4.21) to rule out the occurrence of (i) right away. Indeed, in this case (4.21) and (i) imply that

$$\beta(y) - \beta(z) = n + 1 - 0 > n \geq d(y, z).$$

Since $\operatorname{sgn} \kappa(y) = \operatorname{sgn} \kappa(z) = -1$, this means $\beta$ does not satisfy (SL2), which leads to a contradiction.

On the other hand, suppose (ii) holds. Let $x^{(0)} = z$ [where $z$ satisfies (4.21)], $x^{(r)} = y$ and let $\{x^{(i)}\}_{i=0}^{r-1} \subset \mathcal{S}^n$ be such that $d(x^{(i)}, x^{(i+1)}) = 1$. Note that we can always choose $r \leq n$ because the diameter of $\mathcal{S}^n$ is no greater than $n$. Let $f(i) = \beta(x^{(i)})$; thus $f(0) = 0$ and $f(r) = n$. We say that a "sign change" occurs on the $i$th step if $\operatorname{sgn} \kappa(x^{(i)}) \neq \operatorname{sgn} \kappa(x^{(i+1)})$; we call this sign change positive if $\operatorname{sgn} \kappa(x^{(i)}) < \operatorname{sgn} \kappa(x^{(i+1)})$ and negative otherwise. Since $\kappa(x^{(0)}) < 0$ and $\kappa(x^{(r)}) > 0$ there must be an odd number of sign changes; furthermore, the number of positive sign changes exceeds the number of negative ones by 1. The fact that $\beta$ satisfies property (SL3) implies that $f$ cannot increase on a positive sign change and can increase by at most 2 on a negative sign change. Moreover, since $\beta$ also satisfies (SL2), we know that the value of $f$ can change by at most 1 when there is no sign change. This means that after $r$ steps, $f$ can increase by at most $r - 1$, contradicting the fact that $f(0) = 0$ and $f(r) = n$, as implied by (ii), and $r \leq n$. This establishes the claim (4.20) and hence completes the proof of the lemma. $\square$

We need one more lemma before we can state the main result.

LEMMA 4.7.    *For $n \in \mathbb{N}$, for every $\kappa \in K_n$ and every $\alpha \in A_n(\kappa)$, we have*

$$(4.22) \qquad \langle \kappa, \alpha \rangle \leq \Psi_{n-1}(\kappa').$$

PROOF.    Fix $n \in \mathbb{N}$ and $\kappa \in K_n$. Then for any $\alpha \in A_n(\kappa)$, we prove the relation (4.22) by induction. For $n = 1$, property (SL1) dictates that $0 \leq \alpha(x) \leq \mathbb{1}_{\{\kappa(x) \leq 0\}}$, and so

$$\langle \kappa, \alpha \rangle = \sum_{x \in \mathcal{S} \,:\, \kappa(x) < 0} \kappa(x) \alpha(x) \leq 0 = \Psi_0(\kappa').$$

[Recall that by definition $\Psi_0(\cdot) = 0$.] Now suppose that, for some $j \geq 1$ and $n = j$, (4.22) holds for all $\kappa \in K_n$ and all $\alpha \in A_n(\kappa)$. Pick any $\lambda \in K_{j+1}$ and any $\alpha \in A_{j+1}(\lambda)$ and decompose

$$\langle \lambda, \alpha \rangle = \sum_{y \in \mathcal{S}} \langle \lambda_y, \alpha_y \rangle. \tag{4.23}$$

By property (b) of Lemma 4.4, for all $y \in \mathcal{S}$, $\alpha \in A_{j+1}(\lambda)$ implies $\alpha_y \in B_j(\lambda_y)$. Along with Lemma 4.6, this ensures the existence of $\gamma^{(y)} \in A_j(\lambda_y)$ such that

$$\langle \lambda_y, \alpha_y \rangle \leq \left( \sum_{x \in \mathcal{S}^j} \lambda_y(x) \right)_+ + \langle \lambda_y, \gamma^{(y)} \rangle.$$

Applying the induction hypothesis to $\gamma^{(y)} \in A_j(\lambda_y)$ in the last display and using the trivial identity

$$\sum_{x \in \mathcal{S}^j} \lambda_y(x) = \sum_{u \in \mathcal{S}^{j-1}} \lambda'_y(u),$$

we see that

$$\langle \lambda_y, \alpha_y \rangle \leq \left( \sum_{u \in \mathcal{S}^{j-1}} \lambda'_y(u) \right)_+ + \Psi_{j-1}(\lambda'_y).$$

Together with (4.23) and Lemma 4.3, the last display yields the inequality

$$\langle \lambda, \alpha \rangle \leq \sum_{y \in \mathcal{S}} \left[ \left( \sum_{u \in \mathcal{S}^{j-1}} \lambda'_y(u) \right)_+ + \Psi_{j-1}(\lambda'_y) \right] = \Psi_j(\lambda'),$$

which proves that (4.22) holds for $n = j + 1$. The lemma then follows by induction. $\quad \square$

We now state the main result of this section.

THEOREM 4.8.   *For $n \in \mathbb{N}$, and every $\kappa \in K_n$ and every $\varphi \in \Phi_n$, we have*

$$\langle \kappa, \varphi \rangle \leq \Psi_n(\kappa).$$

PROOF.   Fix $n \geq 1$, $\kappa \in K_n$ and $\varphi \in \Phi_n$. Define the function $\tilde{\varphi}$ on $\mathcal{S}^n$ by

$$\tilde{\varphi}(x) \doteq (\varphi(x) - \mathbb{1}_{\{\kappa(x) > 0\}})_+$$

for $x \in \mathcal{S}^n$. Then

$$\langle \kappa, \varphi \rangle \leq \sum_{x \in \mathcal{S}^n} (\kappa[x])_+ + \langle \kappa, \tilde{\varphi} \rangle$$

holds, since for any $k \in \mathbb{R}$ and $f \in \mathbb{R}_+$,

$$kf \leq (k)_+ + k(f - \mathbb{1}_{\{k>0\}})_+.$$

In addition, $\tilde{\varphi}$ is easily seen to be in $A_n(\kappa)$. An application of Lemma 4.7 then shows that $\langle \kappa, \tilde{\varphi} \rangle \leq \Psi_{n-1}(\kappa')$. Since by definition, $\Psi_n(\kappa) = \sum_{x \in \mathcal{S}^n} (\kappa(x))_+ + \Psi_{n-1}(\kappa')$, we are done.  $\square$

REMARK 4.9.   Recalling the definitions (3.9) and (4.3) of the $\| \cdot \|_{\Phi}$ and $\| \cdot \|_{\Psi}$ norms, respectively, and noting that $\Psi_n(\kappa) \leq \|\kappa\|_{\Psi}$ for all $\kappa \in K_n$, it is clear that Theorem 4.1 is an immediate consequence of Theorem 4.8.

## 5. The martingale bound for finite $\mathcal{S}$.

In Section 3—specifically, Lemma 3.1 and relation (3.7)—we showed that if $\mathcal{S}$ is finite, then given any probability measure $\mathbb{P}$ on $(\mathcal{S}^n, \mathcal{F})$ and a 1-Lipschitz function $\varphi$, for $1 \leq i \leq n-1$,

$$(5.1) \qquad \|V_i(\varphi)\|_{\infty} \leq \max_{y^{i-1} \in \mathcal{S}^{i-1}, w, \hat{w} \in \mathcal{S}} |\langle \kappa[y^{i-1}, w, \hat{w}], \varphi \rangle|,$$

where $V_i(\varphi)$ are the martingale differences defined in (2.1) and the function $\kappa = \kappa[y^{i-1}, w, \hat{w}] \in K_n$ is given explicitly, as in (3.8), by

$$(5.2) \qquad \kappa(x) = \mathbb{1}_{\{x^{i-1}=y^{i-1}\}} (\mathbb{1}_{\{x_i=w\}} \mathbb{P}(x_{i+1}^n \mid y^{i-1}w) - \mathbb{1}_{\{x_i=\hat{w}\}} \mathbb{P}(x_{i+1}^n \mid y^{i-1}\hat{w})).$$

The crux of the proof of Theorem 2.1 for finite $\mathcal{S}$ is the following result, which is proved using Theorem 4.1 of Section 4.

THEOREM 5.1.   *Suppose $\mathcal{S}$ is finite and $\mathbb{P}$ is a probability measure defined on $(\mathcal{S}^n, \mathcal{F})$. Moreover, given any $1 \leq i \leq n-1$, $y^{i-1} \in \mathcal{S}^{i-1}$ and $w, \hat{w} \in \mathcal{S}$, let the function $\kappa[y^{i-1}, w, \hat{w}] \in K_n$ be defined by (5.2) and the coefficients $\eta_{ij}(y^{i-1}, w, \hat{w})$ and $H_{n,i}$ defined by (1.2) and (1.4), respectively. Then for any function $\varphi \in \Phi_n$,*

$$(5.3) \qquad |\langle \kappa[y^{i-1}, w, \hat{w}], \varphi \rangle| \leq 1 + \sum_{j=i+1}^{n} \eta_{ij}(y^{i-1}, w, \hat{w}) \leq H_{n,i}.$$

PROOF.   The second inequality is a direct consequence of the definition of the coefficients $H_{n,i}$. In order to prove the first inequality, first fix $n \in \mathbb{N}$, the measure $\mathbb{P}$ on $(\mathcal{S}^n, \mathcal{F})$, $1 \leq i \leq n-1$, $y^{i-1} \in \mathcal{S}^{i-1}$ and $w, \hat{w} \in \mathcal{S}$, and set $\kappa = \kappa[y^{i-1}, w, \hat{w}]$. Then let $L \doteq n-i+1$ and for $z \in \mathcal{S}^{i-1}$, define the operator $T_z : K_n \to K_L$ as follows: for $\lambda \in K_n$ and $x \in \mathcal{S}^L$,

$$(T_z \lambda)(x) \doteq \lambda(zx).$$

Given $\varphi \in \Phi_n$, note that $T_{y^{i-1}} \varphi \in \Phi_L$ and, due to the structure of $\kappa = \kappa[y^{i-1}, w, \hat{w}]$ given in (5.2), the relation

$$\langle \kappa, \varphi \rangle = \langle T_{y^{i-1}} \kappa, T_{y^{i-1}} \varphi \rangle$$

holds. Combining this with Theorem 4.1 and definitions (3.9) and (4.3) of the norms $\|\cdot\|_\Phi$ and $\|\cdot\|_\Psi$, we observe that

$$|\langle \kappa, \varphi \rangle| \leq \|T_{y^{i-1}\kappa}\|_\Phi \leq \|T_{y^{i-1}\kappa}\|_\Psi = \max_{s \in \{-1,1\}} \psi_L(sT_{y^{i-1}}\kappa).$$

Thus in order to prove the theorem it suffices to show that for $s \in \{-1, 1\}$,

$$(5.4) \qquad \psi_L(s\kappa^{(L)}) \leq 1 + \sum_{j=i+1}^n \eta_{ij}(y^{i-1}, w, \hat{w}),$$

where $\kappa^{(L)} \doteq T_{y^{i-1}}\kappa$.

Now for $\ell = L, L-1, \ldots, 2$, define

$$\kappa^{(\ell-1)} \doteq (\kappa^{(\ell)})',$$

where $(') : K_n \to K_{n-1}$ is the marginal projection operator defined in (4.1). Then $\kappa^{(\ell)} \in K_\ell$ and a direct calculation shows that for $i < j \leq n$ and $x \in \mathcal{S}^{n-j+1}$,

$$(5.5) \quad \kappa^{(n-j+1)}(x) = \mathbb{P}\{X_j^n = x \mid X^i = y^{i-1}w\} - \mathbb{P}\{X_j^n = x \mid X^i = y^{i-1}\hat{w}\}.$$

Since $\kappa^{(n-j+1)}$ is a difference of two probability measures on $\mathcal{S}^{n-j+1}$, we have by (1.18) that

$$\|\kappa^{(n-j+1)}\|_{\text{TV}} = \tfrac{1}{2}\|\kappa^{(n-j+1)}\|_1 = \sum_{x \in \mathcal{S}^{n-j+1}} (\kappa^{(n-j+1)}(x))_+.$$

Together with (5.5), this immediately shows that for $i < j \leq n$,

$$\sum_{x \in \mathcal{S}^{n-j+1}} (\kappa^{(n-j+1)}(x))_+ = \eta_{ij}(y^{i-1}, w, \hat{w}).$$

Now, from the definition of the $\Psi_n$ functional (4.2), we see that

$$\begin{aligned}
\Psi_L(\kappa^{(L)}) &= \sum_{x \in \mathcal{S}^L} (\kappa^{(L)}(x))_+ + \Psi_{L-1}(\kappa^{(L-1)}) \\
&= \sum_{j=1}^L \sum_{x \in S^j} (\kappa^{(j)}(x))_+ \\
&= \sum_{x \in S^L} (\kappa^{(L)}(x))_+ + \sum_{j=i+1}^n \sum_{x \in S^{n-j+1}} (\kappa^{(n-j+1)}(x))_+.
\end{aligned}$$

It follows trivially that $\sum_{x \in S^L} (\kappa^{(L)}(x))_+ \leq 1$. Together, the last three statements show that (5.4) holds when $s = 1$. The inequality (5.4) with $s = -1$ can be established analogously, and hence the proof of the theorem is complete.  □

PROOF OF THEOREM 2.1 FOR FINITE $\mathcal{S}$. Given a finite set $\mathcal{S}$, choose any probability measure $\mathbb{P}$ on $\mathcal{S}^n$. By Remark 2.2, in order to prove Theorem 2.1 it suffices to establish the bound $\|V_i(\varphi)\|_\infty \leq H_{n,i}$ for $1 \leq i \leq n$ only for functions $\varphi \in \Phi_n$. When $i = n$, the bound follows from (3.6) and Lemma 3.1, while for $1 \leq i \leq n - 1$, it can be obtained by taking the maximum of the left-hand side of (5.3) over $y^{i-1} \in \mathcal{S}^{i-1}, w, \hat{w} \in \mathcal{S}$ and combining the resulting inequality with (5.1). $\square$

**6. Extension of the martingale bound to countable $\mathcal{S}$.** In this section we use an approximation argument to extend the proof of Theorem 2.1 from finite $\mathcal{S}$ to countable $\mathcal{S}$. The key to the approximation is the following lemma.

LEMMA 6.1. *Let $\mathcal{S}$ be a countable space and for some $n \in \mathbb{N}$, let $\varphi$ be a 1-Lipschitz function on $\mathcal{S}^n$. Let $\mathbb{P}$ be a probability measure defined on $(\mathcal{S}^n, \mathcal{F})$ such that $\min_{i,\ldots,n} \inf_{y^i \in \mathcal{S}^i \,:\, \mathbb{P}(X^i = y^i) > 0} \mathbb{P}(X^i = y^i) > 0$ as defined in (2.1) and (1.1), respectively. If there exists a sequence of probability measures $\{\mathbb{P}^{(m)}, m \in \mathbb{N}\}$ such that*

$$\lim_{m \to \infty} \|\mathbb{P} - \mathbb{P}^{(m)}\|_{\mathrm{TV}} = 0, \tag{6.1}$$

*then*

$$\lim_{m \to \infty} \bar{\eta}_{ij}^{(m)} = \bar{\eta}_{ij} \quad and \quad \lim_{m \to \infty} \|V_i^{(m)}(\varphi)\|_\infty = \|V_i(\varphi)\|_\infty \tag{6.2}$$

*where, for $m \in \mathbb{N}$, $\{V_i^{(m)}(\varphi)\}$ and $\{\bar{\eta}_{ij}^{(m)}\}$ are the martingale differences and mixing coefficients associated with $\mathbb{P}^{(m)}$, defined in the obvious manner.*

PROOF. The convergence (6.1) automatically implies the convergence in total variation of the conditional distributions $\mathbb{P}^{(m)}(\cdot \mid A)$ to $\mathbb{P}(A)$ for any $A \in \mathcal{S}^n$ with $\mathbb{P}(A) > 0$ (in fact the convergence is uniform with respect to such $A$ under the stipulated condition). As an immediate consequence, we see that $\bar{\eta}_{ij}^{(m)} \to \bar{\eta}_{ij}$ as $m \to \infty$, and (since $\varphi$ is bounded) that $\|V_i^{(m)}(\varphi) - V_i(\varphi)\|_\infty \to 0$, which implies the convergence in (6.2). $\square$

PROOF OF THEOREM 2.1. Suppose $\mathcal{S} = \{s_i : i \in \mathbb{N}\}$ and for $m \in \mathbb{N}$ define $\mathcal{S}_m \doteq \{s_k \in \mathcal{S} : k \leq m\}$. For any probability measure $\mathbb{P}$ on $(\mathcal{S}^n, \mathcal{F})$, define the $m$-truncation of $\mathbb{P}$ to be the following measure on $(\mathcal{S}^n, \mathcal{F})$: for $x \in \mathcal{S}^n$,

$$\mathbb{P}^{(m)}(x) \doteq \mathbb{1}_{\{x \in \mathcal{S}_m^n\}} \mathbb{P}(x) + \mathbb{1}_{\{x = (s_m, s_m, \ldots, s_m)\}} \mathbb{P}(\mathcal{S}^n \setminus \mathcal{S}_m^n). \tag{6.3}$$

Since by construction $\mathbb{P}^{(m)}$ also defines a probability measure on $\mathcal{S}_m^n$ and $\mathcal{S}_m$ is a finite set, it follows from Section 5 [specifically, inequality (5.1) and

Theorem 5.1] that for any 1-Lipschitz function and $1 \le i \le n$,

$$
(6.4) \qquad \|V_i^{(m)}(\varphi)\|_\infty \le 1 + \sum_{j=i+1}^{n} \bar{\eta}_{ij}^{(m)},
$$

where $V_i^{(m)}(\varphi)$ and $\{\bar{\eta}_{ij}^{(m)}\}$ are defined in the obvious fashion.

On the other hand, it is easy to see that the sequence $\mathbb{P}^{(m)}$ converges to $\mathbb{P}$ in total variation norm. Indeed,

$$
(6.5) \qquad \|\mathbb{P} - \mathbb{P}^{(m)}\|_{\mathrm{TV}} \le \sum_{z \in \mathcal{S}^n \setminus \mathcal{S}_m^n} \mathbb{P}(z),
$$

and the r.h.s. must tend to zero, as $m \to \infty$, being the tail of a convergent sum. Theorem 2.1 then follows by taking limits as $m \to \infty$ in (6.4) and applying Lemma 6.1.  $\square$

## 7. Applications of the main result.

7.1. *Bounding $H_{n,i}$ for Markov chains.* Given a countable set $\mathcal{S}$, let $\mathbb{P}$ be a (possibly inhomogeneous) Markov measure on $(\mathcal{S}^n, \mathcal{F})$ with transition kernels $p_k(\cdot \mid \cdot)$, $1 \le k \le n-1$, as defined in (1.6), and let $\theta_i$ be the $i$th contraction coefficient of the Markov chain, as defined in (1.7). The main result of this section is Lemma 7.1, which shows that Theorem 1.2 follows from Theorem 1.1 by establishing the bound (7.3).

For $1 \le k \le n$, let $P^{(k)}$ be the $\mathcal{S} \times \mathcal{S}$ transition probability matrix associated with the $k$th step of the Markov chain: for $1 \le k \le n-1$,

$$
P_{ij}^{(k)} = p_k(j \mid i) \qquad \text{for } i, j \in \mathcal{S}.
$$

Then, using (1.18), the contraction coefficients $\theta_k$ can be rewritten as

$$
(7.1) \qquad \theta_k = \frac{1}{2} \sup_{i,i' \in \mathcal{S}} \sum_{j \in \mathcal{S}} |P_{ij}^{(k)} - P_{i'j}^{(k)}|.
$$

It is a well-known fact that if $\|u\|_1 < \infty$ and $\sum_{i \in \mathcal{S}} u_i = 0$ and $P$ is a transition probability matrix, then

$$
(7.2) \qquad \|u^\top P\|_1 \le \theta_P \|u\|_1,
$$

with $\theta_P$ defined as in the right-hand side of (7.1), but with $P^{(k)}$ replaced by $P$ (for completeness, this fact is included as Lemma A.2 of the Appendix) and where $u^T$ denotes the transpose of $u$.

LEMMA 7.1. *Let $\mathcal{S}$ be a countable set, $\mathbb{P}$ a Markov measure on $(\mathcal{S}^n, \mathcal{F})$ with transition matrices $\{P^{(k)}\}$, $1 \le k \le n-1$, and let $\{\bar{\eta}_{ij}\}$ and $\{\theta_i\}$ be defined by (1.1) and (7.1), respectively. Then for $1 \le i < j \le n$, we have*

$$
(7.3) \qquad \bar{\eta}_{ij} \le \theta_i \theta_{i+1} \cdots \theta_{j-1}
$$

*and so*

$$\|\Delta_n\|_\infty \leq M_n$$

*where $\|\Delta_n\|_\infty$ and $M_n$ are given by (1.3) and (1.8), respectively.*

PROOF.  Let $(X_i)_{1 \leq i \leq n}$ be the coordinate projections on $(\mathcal{S}^n, \mathcal{F}, \mathbb{P})$, that define a Markov chain of length $n$. Fix $1 \leq i \leq n$, $y^{i-1} \in \mathcal{S}^{i-1}$ and $w, \hat{w} \in \mathcal{S}$. Using the relation (1.18) and the definition (1.2) of $\eta_{ij}(y^{i-1}, w, \hat{w})$, we see that

$$\begin{aligned}
&\eta_{ij}(y^{i-1}, w, \hat{w}) \\
&= \|\mathcal{L}(X_j^n \mid X^i = y^{i-1}w) - \mathcal{L}(X_j^n \mid X^i = y^{i-1}\hat{w})\|_{\text{TV}} \\
&= \frac{1}{2} \sum_{x_j^n} |\mathbb{P}\{X_j^n = x_j^n \mid X^i = y^{i-1}w\} - \mathbb{P}\{X_j^n = x_j^n \mid X^i = y^{i-1}\hat{w}\}|.
\end{aligned}$$

However, by the Markov property of $\mathbb{P}$, for any $x_j^n \in \mathcal{S}^{n-j+1}$ and $z \in \mathcal{S}^i$,

$$\mathbb{P}\{X_j^n = x_j^n \mid X^i = z\} = \mathbb{P}\{X_{j+1}^n = x_{j+1}^n \mid X_j = x_j\}\mathbb{P}\{X_j = x_j \mid X_i = z_i\}.$$

Since $\mathbb{P}\{X_{j+1}^n = x_{j+1}^n \mid X_j = x_j\} \leq 1$, we conclude that for $j > i$,

$$\begin{aligned}
&\eta_{ij}(y^{i-1}, w, \hat{w}) \\
&\leq \frac{1}{2} \sum_{x_j \in \mathcal{S}} |\mathbb{P}\{X_j = x_j \mid X_i = w\} - \mathbb{P}\{X_j = x_j \mid X_i = \hat{w}\}| \\
&= \frac{1}{2} \sum_{x_j \in \mathcal{S}} |(e^{(w)} - e^{(\hat{w})})^\top P^{(i)} P^{(i+1)} \cdots P^{(j-1)} e^{(x_j)}| \\
&= \frac{1}{2} \|(e^{(w)} - e^{(\hat{w})})^\top P^{(i)} P^{(i+1)} \cdots P^{(j-1)}\|_1,
\end{aligned}$$

where, for $x \in \mathcal{S}$, $e^{(x)} \in \mathbb{R}^{\mathcal{S}}$ is the unit vector along the $x$ coordinate, that is, for $y \in \mathcal{S}$, $e_x^{(x)} = 1$ and $e_y^{(x)} = 0$ for all $y \in \mathcal{S}, y \neq x$. Since $\|e^{(w)} - e^{(\hat{w})}\|_1 \leq 2$ and the fact that $P^{(k)}$ are transition matrices ensures that $\sum_{i \in \mathcal{S}}[((e^{(w)})^\top P^{(k)})_i - ((e^{(\hat{w})})^\top P^{(k)})_i] = 0$ for all $k \geq 0$, a repeated application of property (7.2) then yields the inequality

$$\eta_{ij}(y^{i-1}, w, \hat{w}) \leq \prod_{k=i}^{j-1} \theta_k.$$

The bound (7.3) follows by taking the supremum of the l.h.s. over all $y^{i-1} \in \mathcal{S}^{i-1}$ and $w, \hat{w} \in \mathcal{S}$. The second bound is a trivial consequence of the first. $\square$

7.2. *Bounding $H_{n,i}$ for hidden Markov chains.* In this section we apply the apparatus developed above to hidden Markov chains. Roughly speaking, the distribution of a hidden Markov chain with finite state space $\mathcal{S}$ (the so-called "observed" state space) is governed by an underlying Markov chain on a "hidden" state space $\hat{\mathcal{S}}$ and a family of stochastic kernels $q_\ell(\cdot|\cdot) \colon \mathcal{S} \times \hat{\mathcal{S}} \mapsto [0,1]$, which quantify the probability of observing the state $x_\ell \in \mathcal{S}$, given that the Markov chain is in the (hidden) state $\hat{s}_\ell$ at the $l$th step. A rigorous definition is as follows. Given transition kernels $p_i(\cdot|\cdot)$, $i = 1, \ldots, n$, on $\hat{\mathcal{S}}$, let $\mu$ be the associated Markov measure on $\hat{\mathcal{S}}^n$: in other words, for $\hat{x} \in \hat{\mathcal{S}}^n$, we have

$$\mu(\hat{x}) = p_0(\hat{x}_1) \prod_{k=1}^{n-1} p_k(\hat{x}_{k+1} \mid \hat{x}_k).$$

Let $\nu$ be the probability measure on $(\hat{\mathcal{S}} \times \mathcal{S})^n$, equipped with the $\sigma$-algebra of all subsets, defined by

$$(7.4) \qquad \nu(\hat{x}, x) = \mu(\hat{x}) \prod_{\ell=1}^{n} q_\ell(x_\ell \mid \hat{x}_\ell),$$

where $q_\ell(\cdot \mid \hat{s})$ is a probability measure on $\mathcal{S}$ for each $\hat{s} \in \hat{\mathcal{S}}$ and $1 \leq \ell \leq n$. It is easy to see that $\nu$ is a Markov measure on $(\hat{\mathcal{S}} \times \mathcal{S})^n$. Indeed, if $Z = ((\hat{X}_i, X_i), 1 \leq i \leq n)$ is a random variable defined on some probability space $(\Omega, \mathcal{F}, \mathbb{P})$ taking values in $(\hat{\mathcal{S}} \times \mathcal{S})^n$ with distribution $\nu$, then the above construction shows that for any $(\hat{x}, x) \in \hat{\mathcal{S}} \times \mathcal{S}$ and $(\hat{y}_1^i, y_1^i) \in (\hat{\mathcal{S}} \times \mathcal{S})^i$,

$$\mathbb{P}\{(\hat{X}_{i+1}, X_{i+1}) = (\hat{x}, x) \mid (\hat{X}_1^i, X_1^i) = (\hat{y}_1^i, y_1^i)\}$$
$$= p_i(\hat{x} \mid \hat{y}_i) q_{i+1}(x \mid \hat{x})$$
$$= \mathbb{P}\{(\hat{X}_{i+1}, X_{i+1}) = (\hat{x}, x) \mid (\hat{X}_i, X_i) = (\hat{y}_i, y_i)\}.$$

The hidden Markov chain measure is then defined to be the $\mathcal{S}^n$-marginal $\rho$ of the distribution $\nu$:

$$(7.5) \qquad \rho(x) = \mathbb{P}\{X = x\} = \sum_{\beta x \in \beta \mathcal{S}^n} \nu(\hat{x}, x).$$

The random process $(X_i)_{1 \leq i \leq n}$ (or measure $\rho$) on $\mathcal{S}^n$ is called a *hidden Markov* chain (resp., measure); it is well known that $(X_i)$ need not be Markov to any order. We will refer to $(\hat{X}_i)$ as the *underlying* chain, which is Markov by construction.

THEOREM 7.1. *Let $(X_i)_{1 \leq i \leq n}$ be a hidden Markov chain, whose underlying chain $(\hat{X}_i)_{1 \leq i \leq n}$ is defined by the transition kernels $p_i(\cdot \mid \cdot)$. Define the* *kth contraction coefficient $\theta_k$ of the underlying chain by*

$$(7.6) \qquad \theta_k = \sup_{\hat{x}, \hat{x}' \in \hat{\mathcal{S}}} \| p_k(\cdot \mid \hat{x}) - p_k(\cdot \mid \hat{x}') \|_{\mathrm{TV}}.$$

*Then the mixing coefficients $\bar{\eta}_{ij}$ associated with the hidden Markov chain $X$ satisfy*

$$(7.7) \qquad \bar{\eta}_{ij} \leq \theta_i \theta_{i+1} \cdots \theta_{j-1},$$

*for $1 \leq i < j \leq n$.*

The proof of Theorem 7.1 is quite straightforward, basically involving a careful bookkeeping of summation indices, rearrangement of sums, and probabilities marginalizing to 1. As in the ordinary Markov case in Section 7.1, the Markov contraction lemma (Lemma A.2) plays a central role.

PROOF OF THEOREM 7.1. For $1 \leq i < j \leq n$, $y_1^{i-1} \in \mathcal{S}^{i-1}$ and $w_i, w_i' \in \mathcal{S}$, by (1.18) we have

$$\eta_{ij}(y_1^{i-1}, w_i, w_i')$$
$$= \tfrac{1}{2} \sum_{x_j^n} |\mathbb{P}\{X_j^n = x_j^n \mid X_1^i = [y_1^{i-1} w_i]\} - \mathbb{P}\{X_j^n = x_j^n \mid X_1^i = [y_1^{i-1} w_i']\}|.$$

Expanding the first conditional probability above, we obtain

$$\mathbb{P}\{X_j^n = x_j^n \mid X_1^i = [y_1^{i-1} w_i]\}$$
$$= \sum_{\hat{s}_j^n} \sum_{\hat{s}_1^i} \mathbb{P}\{X_j^n = x_j^n, (\hat{X}_1^i, \hat{X}_j^n) = (\hat{s}_1^i, \hat{s}_j^n) \mid X_1^i = [y_1^{i-1} w_i]\}$$
$$= \sum_{\hat{s}_j^n} \sum_{\hat{s}_1^i} \frac{\mathbb{P}\{(\hat{X}_1^i, \hat{X}_j^n) = (\hat{s}_1^i, \hat{s}_j^n)\}}{\mathbb{P}\{X_1^i = [y_1^{i-1} w_i]\}}$$
$$\qquad\qquad \times \mathbb{P}\{(X_1^i, X_j^n) = ([y_1^{i-1} w_i], x_j^n) \mid (\hat{X}_1^i, \hat{X}_j^n) = (\hat{s}_1^i, \hat{s}_j^n)\},$$

which can be further simplified using the fact that, by the definition of $\mathbb{P}$,

$$\mathbb{P}\{(X_1^i, X_j^n) = ([y_1^{i-1} w_i], x_j^n) \mid (\hat{X}_1^i, \hat{X}_j^n) = (\hat{s}_1^i, \hat{s}_j^n)\}$$
$$= \nu(x_j^n \mid \hat{s}_j^n) \nu(y_1^{i-1} \mid \hat{s}_1^{i-1}) q_i(w_i \mid \hat{s}_i)$$

and that, due to the Markov property of $X$ and the definition of $\mu$,

$$\mathbb{P}\{(\hat{X}_1^i, \hat{X}_j^n) = (\hat{s}_1^i, \hat{s}_j^n)\} = \mu(\hat{s}_{j+1}^n \mid \hat{s}_j) \mu(\hat{s}_j \mid \hat{s}_i) \mu(\hat{s}_1^i).$$

Expanding $\mathbb{P}\{X_j^n = x_j^n \mid X_1^i = [y_1^{i-1}\, w_i']\}$ in a similar way, we then have

$$
\begin{aligned}
\eta_{ij}&(y_1^{i-1}, w_i, w_i') \\
&= \tfrac{1}{2} \sum_{x_j^n} \left| \sum_{\hat{s}_j^n} \sum_{\hat{s}_1^i} \mu(\hat{s}_{j+1}^n \mid \hat{s}_j) \right. \\
&\qquad\qquad\qquad \left. \times\, \mu(\hat{s}_j \mid \hat{s}_i) \mu(\hat{s}_1^i) \nu(x_j^n \mid \hat{s}_j^n) \nu(y_1^{i-1} \mid \hat{s}_1^{i-1}) \delta(\hat{s}_i) \right|,
\end{aligned}
$$

where we set (recalling that $\rho([y_1^{i-1}\, w]) = \mathbb{P}\{X_1^i = [y_1^{i-1}\, w]\}\ \forall w \in \mathcal{S}$)

$$
\delta(\hat{s}_i) \doteq \frac{q_i(w_i \mid \hat{s}_i)}{\rho([y_1^{i-1}\, w_i])} - \frac{q_i(w_i' \mid \hat{s}_i)}{\rho([y_1^{i-1}\, w_i'])}.
$$

Since $|\sum_{ij} a_i b_j| \le \sum_i a_i |\sum_j b_j|$ for $a_i \ge 0$, $b_i \in \mathbb{R}$, taking the summation over $\hat{s}_j^n$ and the term $\mu(\hat{s}_{j+1}^n|\hat{s}_j)\nu(x_j^n|\hat{s}_j^n)$ outside the absolute value on the right-hand side, interchanging summations over $x_j^n$ and $\hat{s}_j^n$ and then using the fact that $\sum_{\hat{s}_{j+1}^n} \mu(\hat{s}_{j+1}^n|\hat{s}_j) \sum_{x_j^n} \nu(x_j^n|\hat{s}_j^n) = 1$ for every $\hat{s}_j$, we obtain

$$
\begin{aligned}
\eta_{ij}(y_1^{i-1}, w_i, w_i') &\le \tfrac{1}{2} \sum_{\hat{s}_j} \left| \sum_{\hat{s}_1^i} \mu(\hat{s}_1^i)\mu(\hat{s}_j|\hat{s}_i)\nu(y_1^{i-1} \mid \hat{s}_1^{i-1})\delta(\hat{s}_i) \right| \\
&= \tfrac{1}{2} \sum_{\hat{s}_j} \left| \sum_{\hat{s}_i} \mu(\hat{s}_j|\hat{s}_i)\mathbf{h}_{\hat{s}_i} \right|,
\end{aligned}
\tag{7.8}
$$

where $\mathbf{h} \in \mathbb{R}^{\hat{\mathcal{S}}}$ is the vector defined by

$$
\mathbf{h}_{\hat{v}} \doteq \delta(\hat{v}) \sum_{\hat{s}_1^{i-1}} \mu([\hat{s}_1^{i-1}\, \hat{v}])\, \nu(y_1^{i-1} \mid \hat{s}_1^{i-1}).
\tag{7.9}
$$

Let $A^{(i,j)} \in [0,1]^{\hat{\mathcal{S}} \times \hat{\mathcal{S}}}$ be the matrix with $A_{\hat{s},\hat{s}'}^{(i,j)} = \mathbb{P}(\hat{X}_j = \hat{s}' \mid \hat{X}_i = \hat{s}) = \mu(\hat{s}' \mid \hat{s})$ for $s, \hat{s} \in \mathcal{S}$. Then the bound (7.8) can be recast in the form

$$
\eta_{ij}(y_1^{i-1}, w_i, w_i') \le \tfrac{1}{2} \sum_{\hat{s}_j} |(\mathbf{h}^T A^{(i,j)})_{\hat{s}_j}| = \tfrac{1}{2}\|\mathbf{h}^T A^{(i,j)}\|_1.
$$

Since $A^{(i,j)}$ is simply the transition matrix of the Markov chain $\hat{X}$ from step $i$ to $j$, the contraction coefficient of $A^{(i,j)}$ is clearly bounded by $\prod_{k=i}^{j-1} \theta_k$. Therefore, to prove the theorem, it suffices to verify that the assumptions

$$
\sum_{\hat{v} \in \hat{\mathcal{S}}} \mathbf{h}_{\hat{v}} = 0 \quad \text{and} \quad \tfrac{1}{2}\|\mathbf{h}\|_1 \le 1
\tag{7.10}
$$

of the contraction Lemma A.2 are satisfied. Now, expanding (7.9), we have

$$\mathbf{h}_{\hat{v}} = \left( \frac{q_i(w_i \mid \hat{v})}{\rho([y_1^{i-1} w_i])} - \frac{q_i(w_i' \mid \hat{v})}{\rho([y_1^{i-1} w_i'])} \right) \sum_{\hat{s}_1^{i-1}} \mu([\hat{s}_1^{i-1} \hat{v}]) \, \nu(y_1^{i-1} \mid \hat{s}_1^{i-1}).$$

Summing the first term over $\hat{v}$, and using (7.4) and (7.5), we obtain

$$\sum_{\hat{v} \in \hat{\mathcal{S}}} \frac{q_i(w_i \mid \hat{v})}{\rho([y_1^{i-1} w_i])} \sum_{\hat{s}_1^{i-1}} \mu([\hat{s}_1^{i-1} \hat{v}]) \, \nu(y_1^{i-1} \mid \hat{s}_1^{i-1})$$

$$= \frac{1}{\rho([y_1^{i-1} w_i])} \sum_{\hat{s}_1^{i-1}} \nu(y_1^{i-1} \mid \hat{s}_1^{i-1}) \sum_{\hat{v} \in \hat{\mathcal{S}}} q_i(w_i \mid \hat{v}) \mu([\hat{s}_1^{i-1} \hat{v}])$$

$$= \frac{1}{\rho([y_1^{i-1} w_i])} \sum_{\hat{s}_1^{i-1}} \nu(y_1^{i-1} \mid \hat{s}_1^{i-1}) \mathbb{P}\{(\hat{X}_1^{i-1}, X_i) = (\hat{s}_1^{i-1}, w_i)\} = 1.$$

An analogous identity holds for the summation over $\hat{v}$ of the second term, which proves (7.10) and, hence, the theorem. $\square$

Observe that the $\eta$-mixing coefficients of a hidden Markov chain are bounded by the contraction coefficients of the underlying Markov one. One might thus be tempted to pronounce Theorem 7.1 as "obvious" in retrospect, based on the intuition that, conditioned on the hidden Markov chain $\hat{X}$, the observed process $(X_i)_{1 \leq i \leq n}$ is a sequence of independent random variables. Thus, the reasoning might go, all the dependence structure is contained in $\hat{X}_i$, and it is not surprising that the underlying process alone suffices to bound $\bar{\eta}_{ij}$—which, after all, is a measure of the dependence in the process. Such an intuition, however, would be wrong, as it fails to carry over to the case where the underlying process is not Markov. A numerical example of such an occurrence is given in Kontorovich's Ph.D. thesis [15] and, prior to that, in [13], which is also where Theorem 7.1 was first proved. These techniques have been extended further to prove concentration for Markov tree processes; see [14] or [15].

7.3. *Tightness of martingale difference bound for Markov measures.* Given a probability measure $\mathbb{P}$ on $(\mathcal{S}^n, \mathcal{F})$, from (3.2), we know that the associated martingale differences $\{V_i(\varphi)\}$ satisfy

$$(7.11) \qquad \|V_i(\varphi)\|_\infty = \max_{z^i \in \mathcal{S}^i} |V_i(\varphi; z^i)|,$$

where for $1 \leq i \leq n$ and $x^i \in \mathcal{S}^i$,

$$V_i(\varphi; z^i) = \mathbb{E}[\varphi(X) \mid X^i = z^i] - \mathbb{E}[\varphi(X) \mid X^{i-1} = z^{i-1}].$$

Just as $\hat{V}_i(\cdot; y^{i-1}, w, \hat{w})$ could be expressed as an inner product in (3.7) and (3.8), $V_i(\cdot; z^i)$, being a linear functional on $K_n$, also admits a representation in terms of an inner product. Indeed, for $z^i \in \mathcal{S}^i$, we have

$$V_i(\varphi; z^i) = \langle \kappa[z^i], \varphi \rangle$$

where $\kappa = \kappa[z^i] \in K_n$ has the form

$$(7.12) \qquad \kappa(x) = \mathbb{1}_{\{x^i = z^i\}} p(x_{i+1}^n \mid z^i) - \mathbb{1}_{\{x^{i-1} = z^{i-1}\}} p(x_i^n \mid z^{i-1})$$

for $x \in \mathcal{S}^n$. When combined with the definition of the norm $\|\cdot\|_\Phi$ and Theorem 4.1, this shows that

$$\max_{\varphi \in \Phi_n} \|V_i(\varphi)\|_\infty = \max_{z^i \in \mathcal{S}^i} \max_{\varphi \in \Phi_n} |\langle \kappa(z^i), \varphi \rangle| = \max_{z^i \in \mathcal{S}^i} \|\kappa(z^i)\|_\Phi \le \max_{z^i \in \mathcal{S}^i} \|\kappa(z^i)\|_\Psi.$$

It is of interest to ask whether this martingale difference bound is tight, and if so, whether it is possible to obtain a simple description of a class of extremal functions $\varphi$ for which the right-hand side is attained. In this section, we identify such a class when $\mathbb{P}$ is a Markov measure.

The main result is encapsulated in Theorem 7.5, whose statement requires the definition of the BAR class of extremal functions.

DEFINITION 7.2.    A function $\varphi \in \Phi_n$ is said to admit a binary additive representation if there exist functions $\mu_\ell : \mathcal{S} \to \{0, 1\}$, $\ell = 1, \ldots, n$, such that for every $x \in \mathcal{S}^n$,

$$(7.13) \qquad \varphi(x) = \sum_{\ell=1}^n \mu_\ell(x_\ell).$$

In this case, we call $\varphi$ a BAR function and let $\bar{\Phi}_n$ denote the collection of BAR functions in $\Phi_n$.

REMARK 7.3.    Observe that any map $x_1^n \mapsto \mathbb{R}$ of the form (7.13) is 1-Lipschitz and has range in $[0, n]$. Since $\Phi_n$ is an uncountable set while $\bar{\Phi}_n$ is finite, we trivially have $\bar{\Phi}_n \subsetneq \Phi_n$. To get a meaningful size comparison, let us examine the integer-valued members of $\Phi_n$, denoted by $\hat{\Phi}_n \doteq \Phi_n \cap \mathbb{N}^{\mathcal{S}^n}$. For $|\mathcal{S}| \ge 2$, a crude lower bound on the cardinality of $\hat{\Phi}_n$ is

$$|\hat{\Phi}_n| \ge 2^{|\mathcal{S}|^n}.$$

On the other hand, the cardinality of $\bar{\Phi}_n$ is easy to compute exactly:

$$|\bar{\Phi}_n| = 2^{n|\mathcal{S}|}.$$

Thus the vast majority of $\varphi \in \hat{\Phi}_n$ are not BAR functions.

We first begin with a lemma that shows that the norms $\|\cdot\|_\Phi$ and $\|\cdot\|_\Psi$ coincide on the subset of so-called Markov-induced functions in $K_n$. In order to state the lemma, we need to introduce some notation. Fix an (inhomogeneous) Markov measure $\mathbb{P}$ on $(\mathcal{S}^n, \mathcal{F})$ and $z \in \mathcal{S}$, and let $p_0(\cdot)$ and $\{p_k(\cdot\,|\,\cdot) : 1 \le k < n\}$ be the associated initial measure and transition kernels, respectively. In this case, $\kappa = \kappa[z^i]$ in (7.12) can be rewritten as

$$(7.14) \qquad \kappa(x) = \sigma(x_i)\left( \mathbb{1}_{\{x^{i-1} = z^{i-1}\}} \prod_{k=i}^{n-1} p_k(x_{k+1} \,|\, x_k) \right)$$

for $x \in \mathcal{S}^n$, where $\sigma = \sigma[z^i] \in K_1$ is the real-valued function on $\mathcal{S}$ defined by

$$(7.15) \qquad \sigma(y) = \mathbb{1}_{\{y = z_i\}} - p_{i-1}(y \,|\, z_{i-1}).$$

In the case $i = 1$, by our conventions, the above relations reduce to the following:

$$(7.16) \qquad \kappa(x) = \sigma(x_1) \prod_{k=1}^{n-1} p_k(x_{k+1} \,|\, x_k)$$

for $x \in \mathcal{S}^n$, where $\sigma = \sigma(z) \in K_1$ is the function on $\mathcal{S}$ given by

$$\sigma(y) = \mathbb{1}_{\{y = z\}} - p_0(y).$$

For any $n \in \mathbb{N}$ and $\kappa \in K_n$, we say that $\kappa$ is *Markov-induced* if it has the form (7.14), for some collection of transition kernels $\{p_k, 1 \le k < n\}$ with $p_k(z|y) > 0$ for all $1 \le k < n$ and $z, y \in \mathcal{S}$ and function $\sigma \in K_1$.

LEMMA 7.4.    *For any Markov-induced $\kappa \in K_n$, there exists a BAR function $\bar{\varphi} \in \bar{\Phi}_n$ such that*

$$\langle \kappa, \bar{\varphi} \rangle = \Psi_n(\kappa),$$

*and so*

$$|\langle \kappa, \bar{\varphi} \rangle| = \|\kappa\|_\Psi.$$

PROOF.    We shall first prove this result for the case when $i = 1$. In this case, $\kappa$ takes the form (7.16) and satisfies the key property that for $x \in \mathcal{S}^n$,

$$(7.17) \qquad \mathrm{sgn}(\kappa(x)) = \mathrm{sgn}(\sigma(x_1)),$$

meaning that $\mathrm{sgn}(\kappa(x))$ is a function of $x_1$ only. Thus we refer to $\sigma$ as the *sign function* of $\kappa$.

We first claim that for any $\ell \in \mathbb{N}$, if $\kappa^{(\ell)} \in K_\ell$ is of the form (7.16) with some sign function $\sigma^{(\ell)} \in K_1$, then $(\kappa^{(\ell)})' \in K_{\ell-1}$ is Markov-induced with sign function $\sigma^{(\ell-1)}$ given by

$$(7.18) \qquad \sigma^{(\ell-1)}(z) = \sum_{x \in \mathcal{S}} \sigma^{(\ell)}(x) p_1(z \,|\, x) \qquad \text{for } z \in \mathcal{S}.$$

[Here, $(')\colon K_\ell \to K_{\ell-1}$ is the marginal projection operator defined in Section 4.] This is readily verified by observing that for $z \in \mathcal{S}^{\ell-1}$,

$$(\kappa^{(\ell)})'(z) = \sum_{x \in \mathcal{S}} (\kappa^{(\ell)})(xz)$$

$$= \sum_{x \in \mathcal{S}} \sigma^{(\ell)}(x) p_1(z_1 \mid x) \prod_{k=2}^{\ell-1} p_k(z_k \mid z_{k-1})$$

$$= \left( \prod_{k=2}^{\ell-1} p_k(z_k \mid z_{k-1}) \right) \sum_{x \in \mathcal{S}} \sigma^{(\ell)}(x) p_1(z_1 \mid x)$$

$$= \sigma^{(\ell-1)}(z_1) \prod_{k=1}^{\ell-2} p_{k+1}(z_{k+1} \mid z_k),$$

which is of the form (7.16) with sign function $\sigma^{(\ell-1)}$.

Thus, given a Markov-induced $\kappa \in K_n$ with associated sign function $\sigma \in K_1$, first define $\kappa^{(n)} = \kappa$, $\sigma^{(n)} = \sigma$ and, for $\ell = n, \ldots, 2$, let $\kappa^{(\ell-1)} = (\kappa^{(\ell)})'$ and let $\sigma^{(\ell-1)} \in K_1$ be the sign function of $\kappa^{(\ell-1)}$. Then, for $\ell = 1, \ldots, n$, each $\kappa^{(\ell)} \in K_\ell$ satisfies

$$(7.19) \qquad \operatorname{sgn}(\kappa^{(\ell)}(x)) = \operatorname{sgn}(\sigma^{(\ell)}(x_1)).$$

Next, construct the sequence of functions $\mu_1, \ldots, \mu_n$, from the sequence $\kappa^{(1)}$, $\kappa^{(2)}, \ldots, \kappa^{(n)}$ with $\mu_\ell \colon \mathcal{S} \to \{0, 1\}$ given by

$$(7.20) \qquad \mu_\ell(x) = \mathbb{1}_{\{\sigma^{(n-\ell+1)}(x) > 0\}}.$$

Then the function $\bar{\varphi} \colon \mathcal{S}^n \to \mathbb{R}$ defined by

$$(7.21) \qquad \bar{\varphi}(x) = \sum_{\ell=1}^{n} \mu_\ell(x_\ell)$$

for $x \in \mathcal{S}^n$, is easily seen to belong to $\bar{\Phi}_n$. Moreover, note that

$$\langle \kappa, \bar{\varphi} \rangle = \sum_{x \in \mathcal{S}^n} \kappa(x) \bar{\varphi}(x)$$

$$= \sum_{x \in \mathcal{S}^n} \kappa(x) \sum_{\ell=1}^{n} \mu_\ell(x_\ell)$$

$$= \sum_{x \in \mathcal{S}^n} \mu_1(x_1) \kappa(x) + \sum_{x \in \mathcal{S}^n} \mu_2(x_2) \kappa(x) + \cdots + \sum_{x \in \mathcal{S}^n} \mu_n(x_n) \kappa(x)$$

$$= \sum_{x_1^n} \mu_1(x_1) \kappa^{(n)}(x_1^n) + \sum_{x_2^n} \mu_2(x_2) \kappa^{(n-1)}(x_2^n) + \cdots + \sum_{x_n \in \mathcal{S}^1} \mu_n(x_n) \kappa^{(1)}(x_n)$$

$$= \sum_{x_1^n} \mathbb{1}_{\{\sigma^{(n)}(x_1) > 0\}} \kappa^{(n)}(x_1^n) + \sum_{x_2^n} \mathbb{1}_{\{\sigma^{(n-1)}(x_2) > 0\}} \kappa^{(n-1)}(x_2^n) + \cdots$$

$$+ \sum_{x_n} \mathbb{1}_{\{\sigma^{(1)}(x_n) > 0\}} \kappa^{(1)}(x_n)$$

$$= \sum_{x \in \mathcal{S}^n} (\kappa^{(n)}(x))_+ + \sum_{x \in \mathcal{S}^{n-1}} (\kappa^{(n-1)}(x))_+ + \cdots + \sum_{x \in \mathcal{S}} (\kappa^{(1)}(x))_+$$

$$= \Psi_n(\kappa),$$

where the second to last equality uses the property (7.19) and the last equality follows from the definition of the operator $\Psi_n$. This completes the proof of the first statement of the lemma. Due to the definition of the norm $\|\cdot\|_\Psi$, the second statement is a simple consequence of the first.

The case of general $1 \leq i \leq n$ can be dealt with by a corresponding extension of Lemma 7.4 from Markov-induced $\kappa$ of the form (7.16) to $\kappa$ of the form (7.14), which can be achieved by the dimension-reducing technique employed in Section 4. We omit the details.   $\square$

The last lemma immediately implies the following extremal property of BAR functions with respect to martingale differences of Markov measures.

THEOREM 7.5.   *Given a finite set $\mathcal{S}$ and a Markov measure $\mathbb{P}$ with full support on $(\mathcal{S}^n, \mathcal{F})$ as defined in (1.6), for every $1 \leq i \leq n-1$, there exists a BAR function $\bar{\varphi} \in \bar{\Phi}_n$ such that*

$$(7.22) \quad \|V_i(\bar{\varphi})\|_\infty = \max_{\varphi \in \Phi_n} \|V_i(\varphi)\|_\infty = \max_{z^i \in \mathcal{S}^i} \|\kappa[z^i]\|_\Phi = \max_{z^i \in \mathcal{S}^i} \|\kappa[z^i]\|_\Psi.$$

PROOF.   Given any $\varphi \in \Phi_n$, let $z^i$ be the element of $\mathcal{S}^i$ that achieves the maximum in the right-hand side of (7.11). Then the discussion at the beginning of the section, along with Lemma 7.4 and Theorem 4.1, shows that (when $i = 1$) there exists a BAR function $\bar{\varphi}$ such that

$$|V_i(\bar{\varphi}; z^i)| = |\langle \kappa[z^i], \bar{\varphi} \rangle| = \|\kappa[z^i]\|_\Psi \geq \|\kappa[z^i]\|_\Phi \geq |\langle \kappa[z^i], \varphi \rangle| = \|V_i(\varphi)\|_\infty.$$

Taking the maximum over $z^i \in \mathcal{S}^i$, we conclude that

$$\|V_i(\bar{\varphi})\|_\infty = \max_{z^i \in \mathcal{S}^i} \|\kappa[z^i]\|_\Psi \geq \max_{z^i \in \mathcal{S}^i} \|\kappa[z^i]\|_\Phi \geq \|V_i(\varphi)\|_\infty.$$

Taking the maximum of the left-hand side over BAR functions $\bar{\varphi}$ (and, without loss of generality, denoting a maximizing function there again by $\bar{\varphi}$), and then taking the maximum over the right-hand side over functions $\varphi \in \Phi_n$, the fact that $\bar{\Phi}_n \subset \Phi_n$ shows that the inequalities can be replaced by equalities and hence (7.22) follows.   $\square$

## APPENDIX

**A.1. The norms $\|\cdot\|_\Phi$ and $\|\cdot\|_\Psi$.** The following result, while not directly used in the paper, may be of independent interest.

LEMMA A.1. *For $n \geq 1$, both the functionals $\|\cdot\|_\Phi$ and $\|\cdot\|_\Psi$ described in* (3.9) *and* (4.3), *respectively, define norms on $K_n$.*

PROOF. It follows trivially from the definition

$$\|\kappa\|_\Phi = \max_{\varphi \in \Phi_n} |\langle \kappa, \varphi \rangle| = \max_{\varphi \in \Phi_n} \left| \sum_{x \in \mathcal{S}} \kappa(x)\varphi(x) \right|$$

that for $\kappa \in K_n$, $\|\kappa\|_\Phi \geq 0$, $\|\kappa\|_\Phi = 0$ if and only if $\kappa \equiv 0$ (to see this, choose $\varphi(x) = \mathbb{1}_{\{\kappa(x)>0\}}$) and $\|a\kappa\|_\Phi = |a|\|\kappa\|_\Phi$ for $a \in \mathbb{R}$. Last, the triangle inequality $\|\kappa_1 + \kappa_2\|_\Phi \leq \|\kappa_1\|_\Phi + \|\kappa_2\|_\Phi$ follows immediately from the linearity of $\langle \cdot, \cdot \rangle$ in the first variable and the fact that $|\cdot|$ satisfies the triangle inequality. This shows that $\|\cdot\|_\Phi$ defines a norm on $K_n$.

We now consider the functional

$$\|\kappa\|_\Psi = \max\{\Psi_n(\kappa), \Psi_n(-\kappa)\}$$

with the operator $\Psi_n$ defined recursively through the relation

$$\Psi_n(\kappa) \doteq \sum_{x \in \mathcal{S}^n} (\kappa(x))_+ + \Psi_{n-1}(\kappa'),$$

with $\kappa' \in K_{n-1}$ given by $\kappa'(y) = \sum_{x_1 \in \mathcal{S}} \kappa(x_1 y)$. The fact that $\Psi_0 \equiv 0$, along with the above recursion relation, immediately guarantees that for all $\ell = 0, 1, \ldots,$ and $\lambda \in K_\ell$, $\Psi_\ell(\lambda) \geq 0$ and $\Psi_\ell(\lambda) \geq \sum_{x \in \mathcal{S}} (\lambda(x))_+$. If $\kappa(x) \neq 0$ for some $x \in \mathcal{S}$, then the latter quantity is strictly positive for either $\lambda = \kappa$ or $\lambda = -\kappa$, which implies that $\|\kappa\|_\Psi = 0$ if and only if $\kappa \equiv 0$. The homogeneity property of the norm $\|\cdot\|_\Psi$ follows from the corresponding property, for the operator $\Psi_n$—namely, $\Psi_n(a\kappa) = a\Psi_n(\kappa)$ for $a > 0$. Last, the triangle inequality is a consequence of the property $\Psi_n(\kappa_1 + \kappa_2) \leq \Psi_n(\kappa_1) + \Psi_n(\kappa_2)$ for every $\kappa_1, \kappa_2 \in K_n$, which can be deduced using the subadditivity of the function $f(z) = (z)_+$, the fact that $\Psi_0$ trivially satisfies the triangle inequality and induction. $\square$

**A.2. Contraction lemma.** For completeness, we include the elementary proof of a bound that was used in the proof of Lemma 7.1. For finite $\mathcal{S}$, a proof of this result goes back to Markov [18] (see Section 5 of that work, or Lemma 10.6(ii) of [3]). We recall that $u^\top$ denotes the transpose of the vector $u \in \mathbb{R}^\mathcal{S}$.

Lemma A.2.   *Let $\mathcal{S}$ be a countable set, $u \in \mathbb{R}^{\mathcal{S}}$ be such that $\sum_{i \in \mathcal{S}} u_i = 0$ and $\|u\|_1 < \infty$, and let $P$ be an $\mathcal{S} \times \mathcal{S}$ matrix such that $u^{\top}P$ is well defined. Then*

$$
\|u^{\top}P\|_1 \leq \theta_P \|u\|_1, \tag{A.1}
$$

*where $\theta_P$ is the contraction coefficient of $P$:*

$$
\theta_P = \tfrac{1}{2} \sup_{i,i' \in \mathcal{S}} \sum_{j \in \mathcal{S}} |P_{ij} - P_{i'j}|. \tag{A.2}
$$

Proof.   Let $y_i = |u_i|$, and define $I_+$, $I_-$ as follows: $I_+ = \{i \in \mathcal{S} : u_i > 0\}$ and $I_- = \{i \in \mathcal{S} : u_i < 0\}$. Then for any finite $J \subseteq \mathcal{S}$,

$$
\sum_{j \in J} (u^{\top}P)_j = \sum_{j \in J} \left[ \sum_{i \in I_+} y_i P_{ij} - \sum_{i \in I_-} y_i P_{ij} \right]
$$
$$
= \sum_{i \in I_+} Q_i y_i - \sum_{i \in I_-} Q_i y_i
$$

where $Q_i = \sum_{j \in J} P_{ij}$. Thus, we obtain

$$
\left| \sum_{j \in J} (u^{\top}P)_j \right| = \left| \sum_{i \in I_+} Q_i y_i - \sum_{j \in I_-} Q_i y_i \right|
$$
$$
\leq \left| \left( \sup_{k \in I_+} Q_k \right) \sum_{i \in I_+} y_i - \left( \inf_{k \in I_-} Q_k \right) \sum_{i \in I_-} y_i \right|
$$
$$
= \tfrac{1}{2} \|u\|_1 \left| \sup_{k \in I_+} Q_k - \inf_{k \in I_-} Q_k \right|
$$
$$
\leq \tfrac{1}{2} \|u\|_1 \sup_{i,i'} \left| \sum_{j \in J} (P_{ij} - P_{i'j}) \right|
$$
$$
\leq \theta_P \|u\|_1.
$$

Taking the supremum of the l.h.s. over all finite $J \subseteq \mathcal{S}$ yields the result.   □

**Acknowledgment.**   The authors would like to thank John Lafferty for his support during the course of this work.

Dᴇᴘᴀʀᴛᴍᴇɴᴛ ᴏғ Mᴀᴛʜᴇᴍᴀᴛɪᴄs
Wᴇɪᴢᴍᴀɴɴ Iɴsᴛɪᴛᴜᴛᴇ ᴏғ Sᴄɪᴇɴᴄᴇ
Rᴇʜᴏᴠᴏᴛ 76100
Isʀᴀᴇʟ
E-ᴍᴀɪʟ: aryehk.kontorovich@weizmann.ac.il

Dᴇᴘᴀʀᴛᴍᴇɴᴛ ᴏғ Mᴀᴛʜᴇᴍᴀᴛɪᴄᴀʟ Sᴄɪᴇɴᴄᴇs
Cᴀʀɴᴇɢɪᴇ Mᴇʟʟᴏɴ Uɴɪᴠᴇʀsɪᴛʏ
Pɪᴛᴛsʙᴜʀɢʜ, Pᴇɴɴsʏʟᴠᴀɴɪᴀ 15213
USA
E-ᴍᴀɪʟ: kramanan@math.cmu.edu